\nonstopmode
\documentclass[letterpaper, 10pt]{amsart}
\usepackage{latexsym}
\usepackage{fancyhdr}
\usepackage{amsmath, amssymb}
\usepackage[ansinew]{inputenc}
\usepackage[all]{xy}
\usepackage{verbatim}
\usepackage{pdflscape}
\usepackage{longtable}
\usepackage{rotating}
\usepackage{psfrag}
\usepackage{epsfig}
\usepackage{pdfsync}
\usepackage{hyperref}
\usepackage{subfigure}
\swapnumbers
\theoremstyle{plain}
\newtheorem*{lemma*}{Lemma}
\newtheorem{lemma}[subsection]{Lemma}
\newtheorem*{theorem*}{Theorem}
\newtheorem{theorem}[subsection]{Theorem}
\newtheorem*{proposition*}{Proposition}
\newtheorem{proposition}[subsection]{Proposition}
\newtheorem*{corollary*}{Corollary}
\newtheorem{corollary}[subsection]{Corollary}
\newtheorem*{claim*}{Claim}

\theoremstyle{definition}
\newtheorem*{definition*}{Definition}

\newtheorem*{example*}{Example}
\newtheorem{example}[subsection]{Example}
\newtheorem*{algorithm*}{Algorithm}
\newtheorem*{remark*}{Remark}
\newtheorem*{remarks*}{Remarks}

\newtheorem{remark}[subsection]{Remark}
\newtheorem{remarks}[subsection]{Remarks}

\numberwithin{equation}{subsection}
\newenvironment{demo}[1]{\par\smallskip\noindent{\bf #1.}}{\par\smallskip}

\newcommand{\East}[2]{-\raisebox{0.1pt}{$\mkern-16mu\frac{\;\;#1\;}{\;\;#2\;}\mkern-16mu$}\to}

\def\al{\alpha}
\def\be{\beta}
\def\ga{\gamma}
\def\de{\delta}
\def\ep{\epsilon}

\def\la{\lambda}

\def\rh{\rho}

\def\si{\sigma}

\def\ta{\tau}

\def\ph{\phi}
\def\vh{\varphi}
\def\ch{\chi}
\def\ps{\psi}
\def\om{\omega}

\def\De{\Delta}

\def\La{\Lambda}

\def\Ph{\Phi}
\def\Ps{\Psi}
\def\Om{\Omega}

\def\C{\mathbb{C}}

\def\I{\mathbb{I}}

\def\N{\mathbb{N}}

\def\R{\mathbb{R}}

\def\cC{\mathcal{C}}

\def\cF{\mathcal{F}}

\def\cH{\mathcal{H}}
\def\cI{\mathcal{I}}

\def\cL{\mathcal{L}}

\def\cO{\mathcal{O}}

\def\cW{\mathcal{W}}

\def\p{\partial}

\renewcommand{\Re}{\mathrm{Re}}
\renewcommand{\Im}{\mathrm{Im}}
\def\<{\langle}
\def\>{\rangle}
\renewcommand{\o}{\circ}

\let\on=\operatorname

\title[Perturbation of polynomials and normal matrices]
{Quasianalytic multiparameter perturbation of polynomials and normal matrices}

\author[A.\ Rainer]
{Armin Rainer}

\address{Armin Rainer: Department of Mathematics, University of Toronto, 
40 St.\ George Street, Toronto, Ontario, Canada M5S 2E4}

\email{armin.rainer@univie.ac.at}

\begin{document}

\begin{abstract}
We study the regularity of the roots of 
multiparameter families of complex univariate monic polynomials 
$
P(x)(z) = z^n + \sum_{j=1}^n (-1)^j a_j(x) z^{n-j}
$ 
with fixed degree $n$ whose coefficients belong to a certain subring $\cC$ of $C^\infty$-functions.
We require that $\cC$ includes polynomial but excludes flat functions (quasianalyticity) 
and is closed under composition, derivation, division by a coordinate, and taking the inverse.
Examples are quasianalytic Denjoy--Carleman classes, in particular, the class of real analytic functions $C^\om$.

We show that there exists a locally finite covering $\{\pi_k\}$ of the parameter space, where each $\pi_k$ is a 
composite of finitely many $\cC$-mappings each of which is either a local blow-up with smooth center or a local power substitution
(in coordinates given by $x \mapsto (\pm x_1^{\ga_1},\ldots,\pm x_q^{\ga_q})$, $\ga_i \in \N_{>0}$),
such that, for each $k$, 
the family of polynomials $P \o \pi_k$ admits a $\cC$-parameterization of its roots. 
If $P$ is hyperbolic (all roots real), then local blow-ups suffice.  

Using this desingularization result, we prove that the roots of $P$ can be parameterized by 
$SBV_{\on{loc}}$-functions 
whose classical gradients exist almost everywhere and belong to $L^1_{\on{loc}}$.
In general the roots cannot have gradients 
in $L^p_{\on{loc}}$ for any $1 < p  \le \infty$. Neither can the roots be in $W_{\on{loc}}^{1,1}$ or $VMO$.

We obtain the same regularity properties for the eigenvalues and the eigenvectors of $\cC$-families of normal matrices.
A further consequence is that every continuous subanalytic function belongs to $SBV_{\on{loc}}$.
\end{abstract}

{\noindent{\small\rm To appear in Trans.\ Amer.\ Math.\ Soc.}} 

\thanks{The author was supported by the Austrian Science Fund (FWF), Grant J2771}
%\dedicatory{dedicatory}
\date{January 28, 2010}

\maketitle

\section{Introduction}

Let us consider a family of univariate monic polynomials
\[
P(x)(z) = z^n + \sum_{j=1}^n (-1)^j a_j(x) z^{n-j}
\]
where the coefficients $a_j : U \to \C$ (for $1 \le j \le n$) are complex
valued functions defined in an open
subset $U \subseteq \R^q$. 
If the coefficients $a_j$ are regular (of some kind)
it is natural to ask whether the roots of $P$ can be arranged regularly
as well, i.e.,
whether it is possible to find $n$ regular functions $\la_j : U \to \C$ (for $1
\le j \le n$) such that $\la_1(x),\ldots,\la_n(x)$ represent the roots of
$P(x)(z)=0$ for each $x \in U$.

This perturbation problem has been intensively studied under the following
additional assumptions:
\begin{enumerate}
\item The parameter space is one dimensional: $q=1$. 
\item The polynomials $P(x)$ are hyperbolic, i.e., all roots of $P(x)$ are
real.
\end{enumerate}

If both of these conditions are satisfied, there exist real analytic
parameterizations of
the roots of $P$ if its coefficients $a_j$ are real analytic, by a classical
theorem due to Rellich \cite{Rellich37I}.
If all $a_j$ are smooth
($C^\infty$) and no two of the
increasingly ordered (hence) continuous roots meet of infinite order of flatness,
then there exist smooth parameterizations of the roots, by \cite{AKLM98}. 
Without additional condition we cannot hope for smooth roots.
By \cite{RainerOmin}, smooth roots exist if the coefficients are smooth and definable in some o-minimal expansion of the real field, 
which implies that not flat contact but oscillatory behavior is responsible for the loss of smoothness.
The roots may always be chosen $C^1$ (resp.\ twice differentiable)
provided that the $a_j$ are in $C^{2n}$ (resp.\ $C^{3n}$), see \cite{Mandai85} and \cite{KLM04}. 
Recently, the assumptions in this statement have been refined to $C^n$ (resp.\ $C^{2n}$) by \cite{ColombiniOrruPernazza08}.
It is then best possible in both hypothesis and conclusion as shown by examples (e.g.\ in \cite{ColombiniOrruPernazza08} and \cite{BBCP06}).
Sharp sufficient conditions, in terms of the differentiability of the coefficients and the order of contact of the roots, for the existence of $C^p$-roots
($p \in \N$)  
are found in \cite{RainerOmin}.

If the polynomials $P(x)$ are hyperbolic and all $a_j$ are in $C^n$, 
but the parameter space is multidimensional ($q >1$), then the roots of $P$
may still be parameterized by locally Lipschitz functions (by ordering them increasingly for instance). 
This follows from the
fundamental results of Bronshtein \cite{Bronshtein79} and (alternatively)
Wakabayashi \cite{Wakabayashi86} (which also constitute the main part in the proof of all but the last of the finite differentiability statements above). 
For a detailed presentation of those see \cite{RainerMONO}.
A different and easier proof for the
partial case that the
coefficients $a_j$ are real analytic was recently given by Kurdyka and
Paunescu \cite{KurdykaPaunescu08}. 
In that paper the real analytic multiparameter perturbation theory of hyperbolic polynomials $P$ and symmetric matrices $A$ is studied. 
It is shown that there exists a modification $\Ph : W \to U$, namely a locally finite composition
of blow-ups with smooth centers, such that the roots of $P \o \Ph$ can be
locally parameterized by real analytic functions, and $A \o \Ph$ is real analytically diagonalizable. 
For further results on the perturbation problem of hyperbolic polynomials see (among others) \cite{Glaeser63R},
\cite{Dieudonne70}, \cite{CC04}, and \cite{LR07}.

The one parameter case $q=1$, but with the hyperbolicity assumption dropped, was
treated in \cite{RainerAC}.
In that case continuous parameterizations of
the roots still exist given that
the coefficients $a_j$ are continuous (e.g.\ Kato \cite[II 5.2]{Kato76}).
If all $a_j$ are smooth and no two of
the continuously chosen
roots meet of infinite order of flatness, then any continuous
parameterization of the roots is locally
absolutely continuous.  
Absolute continuity is the best one can expect, see \ref{optimal}.
This theorem follows from the (Puiseux type) proposition that for any $x_0$ there exists
an integer $N$ such that
$x \mapsto P(x_0\pm(x-x_0)^N)$ admits smooth parameterizations of its roots
near $x_0$.
It seems unknown whether the roots still can be arrange locally absolutely continuously if the condition on the order of contact is omitted.
Spagnolo \cite{Spagnolo99} gave an affirmative answer for degree 2 and 3 polynomials
(degree  4 is announced).

In the present paper we study smooth multiparameter perturbations
of complex polynomials, i.e.,
without the restrictions $(1)$ and $(2)$. 
It is easy to see that every choice of the roots of a bounded family $P$ of polynomials is bounded as well (proposition \ref{bdroot}).
By a theorem due to
Ostrowski \cite{Ostrowski40}, 
for a continuous family $P$ of polynomials, the set of all roots still 
is continuous and satisfies a H\"older condition of order $1/n$.
But in general there may not exist continuous parameterizations of 
the single roots as in the one dimensional or
hyperbolic case. For instance, $P(x_1,x_2)(z) = z^2-(x_1+i x_2)$,
with $x_1,x_2 \in \R$ and $i = \sqrt{-1}$.
Nevertheless, the roots of $P$ may have some other regularity properties.

We show the following (theorem \ref{cCperturb}): 
Let $\cC$ be a certain class of $C^\infty$-functions (specified below).
If the coefficients $a_j$ of $P$ are $\cC$-functions on a $\cC$-manifold $M$, 
then for each compact subset $K \subseteq M$ there exist:
\begin{enumerate}
\item[(a)] a neighborhood $W$ of $K$, and
\item[(b)] a finite covering $\{\pi_k : U_k \to W\}$ of $W$ by $\cC$-mappings, 
where each $\pi_k$ is a composite of finitely many mappings each of which is either a local blow-up 
with smooth center or a local power substitution,
\end{enumerate}
such that, for all $k$, 
the family of polynomials $P \o \pi_k$ allows a $\cC$-parameterization of its roots on $U_k$. 
If $P$ is hyperbolic, then local blow-ups suffice (theorem \ref{hypcCperturbth}).
A local blow-up over an open subset $U \subseteq M$ is a blow-up over $U$ composed with the inclusion of $U$ in $M$. 
A local power substitution is the composite of the inclusion of a coordinate chart $W$ in $M$ and a mapping $V \to W$ given 
in local coordinates by
\[
(x_1,\ldots,x_q) \mapsto ((-1)^{\ep_1} x_1^{\ga_1},\ldots,(-1)^{\ep_q}x_q^{\ga_q})
\]
for some $\ga \in (\N_{>0})^q$ and all $\ep \in \{0,1\}^q$.
(See \ref{not} for a precise explanation of these notions.)

The proof uses resolution of singularities.  
Accordingly, $\cC$ is a class of $C^\infty$-functions admitting resolution of singularities. Due to Bierstone and Milman 
\cite{BM04} (and \cite{BM97}), it suffices that $\cC$ is a subring of $C^\infty$ that includes polynomial 
but excludes flat functions (quasianalyticity) and is closed under composition, differentiation, division by a coordinate, and taking the inverse 
(see section \ref{secC}). 
For instance, $\cC$ may be any quasianalytic Denjoy--Carleman class $C^M$, 
where the weight sequence $M$ satisfies some mild conditions (see section \ref{secDC}).
In particular, $\cC$ can be the class of real analytic functions $C^\om$. 
Hence, in the hyperbolic case, we recover a version of the aforementioned theorem 
due to Kurdyka and Paunescu \cite{KurdykaPaunescu08}. 

The above result (theorem \ref{cCperturb}) enables us to investigate the regularity of 
the roots of $\cC$-families of polynomials $P$. 
We show:
\begin{enumerate}
\item[(i)] The roots of $P$ allow a parameterization by ``piecewise Sobolev $W^{1,1}_{\on{loc}}$'' functions. More precisely,
the roots of $P$ can locally be chosen of class $\cC$ outside of a closed nullset of finite $(q-1)$-dimensional Hausdorff measure 
such that its classical gradient belongs to $L^1_{\on{loc}}$ (theorem \ref{Wloc}).
\item[(ii)] The roots of $P$ allow a parameterization in $SBV_{\on{loc}}$ (theorem \ref{BVth}). 
\end{enumerate}
Note that (i) implies (ii) (see section \ref{secBV}).
Simple examples show that the conclusion in (i) is best possible:
In general we cannot expect that the roots of $P$ admit arrangements having gradients
in $L^p_{\on{loc}}$ for any $1 < p \le \infty$ (see \ref{optimal}). 
In contrast to the one parameter case (see \cite{RainerAC} and \ref{1dim}), 
multiparameter families of polynomials do in general not allow roots in $W^{1,1}_{\on{loc}}$ 
(see the polynomial counter-example in \ref{sqroot}) or in $VMO$ (see \ref{VMO}).

The question for optimal assumptions is open. For instance,
it is unknown whether (ii) still holds when the coefficients of $P$ are just $C^\infty$-functions. 
That problem requires different methods.

Table \ref{summary} on page \pageref{summary} provides a summary of the most important results on the perturbation theory of polynomials. 

In section \ref{secmatrix}  
we deduce consequences for the perturbation theory of normal matrices.
There will be applications to the perturbation theory of unbounded normal operators with compact resolvents and
common domain of definition as well. It requires a differential calculus for quasianalytic classes beyond Banach spaces 
(see \cite{KMRc} for the case of non-quasianalytic Denjoy--Carleman classes). 
This will be taken up elsewhere (see \cite{KMRq} and \cite{KMRp}).
Our results generalize theorems obtained in \cite{KurdykaPaunescu08} and \cite{RainerAC}.
For more on the perturbation theory of linear operators
consider Rellich \cite{Rellich37I, Rellich37II, Rellich39III, Rellich40IV, Rellich42V, Rellich69},
Kato \cite{Kato76}, Baumg\"artel \cite{Baumgaertel72}, and also \cite{AKLM98}, \cite{KM03}, and \cite{KMR}.

We prove the following (theorem \ref{qamatrix}):
Let $A = (A_{ij})_{1 \le i,j \le n}$ be a family of normal complex matrices, where the entries $A_{ij}$ are $\cC$-functions on a 
$\cC$-manifold $M$. 
Then, for each compact subset $K \subseteq M$, there exist a neighborhood $W$ of $K$ and a finite covering 
$\{\pi_k : U_k \to W\}$ of $W$ of the type described in (b), 
such that, for all $k$, the family of normal matrices
$A \o \pi_k$ allows $\cC$-parameterizations of its eigenvalues and its eigenvectors.
If $A$ is a family of Hermitian matrices, then local blow-ups suffice.
Both a nonflatness condition (such as quasianalyticity) and normality of the matrices $A(x)$ are necessary for the desingularization of the eigenvectors
(see \ref{flat} and \ref{normal}).

We conclude that the eigenvalues and the eigenvectors of a $\cC$-family of normal complex matrices $A$ locally admit parameterizations 
by ``piecewise Sobolev $W^{1,1}_{\on{loc}}$'' functions (in the sense of (i)) and, thus, by 
$SBV_{\on{loc}}$-functions (theorem \ref{regmatrix}).

A further application of the method developed in this paper is given 
in section \ref{secsuban}: Any continuous subanalytic function belongs to $SBV_{\on{loc}}$.

\subsection*{Notation}
We use $\N = \N_{>0} \cup \{0\}$. 
Let $\al=(\al_1,\ldots,\al_q) \in \N^q$ and $x = (x_1,\ldots,x_q) \in \R^q$.
We write $\al!=\al_1! \cdots \al_q!$, $|\al|= \al_1 +\cdots+ \al_q$, $x^\al = x_1^{\al_1}\cdots x_q^{\al_q}$, and 
$\p^\al=\p^{|\al|}/\p x_1^{\al_1} \cdots \p x_q^{\al_q}$. 
We shall also use $\p_i = \p/\p x_i$. If $\al,\be \in \N^q$, then $\al \le \be$ means $\al_i \le \be_i$ for all $1 \le i \le q$.

Let $U \subseteq \R^q$ be an open subset.
For a function $f \in C^\infty(U)$ we denote by 
$\hat f_a \in \cF_q$ its Taylor series at $a \in U$, i.e., 
\[
\hat f_a(x) = \sum_{\al \in \N^q} \frac{1}{\al!} \p^\al f(a) x^\al,
\]
where $\cF_q$ denotes the ring of power series in $q$ variables. %(with real or complex coefficients).

$\on{S}_n$ denotes the symmetric group on $\{1,2,\ldots,n\}$.

We denote by $\cH^q$ (resp.\ $\cL^q$) the $q$-dimensional Hausdorff (resp.\ Lebesgue) measure. 
We also use $|X| = \cL^q(X)$ and $\int_X f(x) dx = \int_X f(x) d \cL^q(x)$.
We write $\mathbf{1}_X$ for the indicator function of a set $X$. 
For $x \in \R^q$, $B_r(x) = \{y \in \R^q : |x-y|<r\}$ is the open ball with center $x$ and radius $r$ with respect to the Euclidean metric.  

All manifolds in this paper are assumed to be Hausdorff, paracompact, and finite dimensional.

{\scriptsize

\begin{landscape}
\setlength{\LTcapwidth}{6in}
\begin{longtable}{|l|c|l|c|l|l|}
\caption{\label{summary}Let
$
P(x)(z) = z^n + \sum_{j=1}^n (-1)^j a_j(x) z^{n-j}
$
be a family of polynomials with coefficients $a_j : \R^q \to \C$ (for $1 \le j \le n$). 
The table provides a (by no means exhaustive) summary of the most important results 
concerning the existence of parameterizations of the roots of $P$ of some regularity, given that $P$ fulfills certain conditions.
The regularity of the roots is in general best possible under the respective conditions on $P$, which might partly not be optimal. 
`Definable' refers to an arbitrary but fixed o-minimal expansion of the real field.
By $\cC$ we mean a class of $C^\infty$ functions satisfying (\ref{cC}.1)--(\ref{cC}.6). For a definition of $\cW^\cC$ see \ref{classW}. 
Normal nonflatness is introduced in \ref{1dim}.
And $s$ is maximal with the property that $\tilde \De_s(P)\ne 0$, where $\tilde \De_s(P)$ is given by \eqref{eqdel}.} \\
\endfirsthead
\hline\endhead
\hline 
& & & & & \\ [-.7ex]
{\bf Extra conditions} & {\&} & {\bf Coefficients} & $\implies$ & {\bf Roots} & {\bf Reference} \\ [0.5ex]
\hline
\hline
& & & & & \\ [-1.5ex]
$q=1$ &  & continuous & & continuous & \cite[II 5.2]{Kato76}\\
[0.5ex]\hline 
& & & & & \\ [-1.5ex]
$q=1$ &  & continuous \& definable & & $AC_{\on{loc}}$ \& definable & \cite{RainerOmin}\\
[0.5ex]\hline 
& & & & & \\ [-1.5ex]
$q=1$ &  & $C^\infty$ \& normally nonflat & & local desingularization by $x \mapsto \pm x^\ga$ ($\ga \in \N_{>0}$), & \cite{RainerAC}\\
& & & & $AC_{\on{loc}}$ \& no two distinct roots meet $\infty$-flat & \\
[0.5ex]\hline 
& & & & & \\ [-1.5ex]
$q=1$ \& $n=2,3,(4)$ &  & $C^\infty$ & & $AC_{\on{loc}}$ & \cite{Spagnolo99}\\
[0.5ex]\hline \hline
& & & & & \\ [-1.5ex] 
& & bounded & & bounded  & proposition \ref{bdroot} \\  
[0.5ex]\hline
& & & & & \\ [-1.5ex]
& & continuous & & continuous as a set, &\cite{Ostrowski40} \\
& & & & fulfill a H\"older condition of order $1/n$& \\
[0.5ex]\hline
& & & & & \\ [-1.5ex]
& & $\cC$ & & local desingularization by finitely many & theorem \ref{cCperturb}\\
& & (resp.\ continuous \& & & local blow-ups with smooth center and & (resp.\ theorem \ref{suban}) \\
& & subanalytic) & & local power substitutions (in the sense of \ref{not}), & \\
& & & & $\cW^{\cC}_{\on{loc}}$ \& $SBV_{\on{loc}}$ & theorems \ref{Wloc} \& \ref{BVth} \\
[0.5ex]\hline
& & & & & \\ [-1.5ex]
hyperbolic \& $q=1$ &  & $C^\om$ (resp.\ $\cC$) & & $C^\om$ (resp.\ $\cC$) & \cite{Rellich37I} (resp.\ corollary \ref{hypcCperturbcor}) \\
[0.5ex]\hline 
& & & & & \\ [-1.5ex]
hyperbolic \& $q=1$ &  & $C^\infty$ \& normally nonflat & & $C^\infty$ \& no two distinct roots meet $\infty$-flat & \cite{AKLM98}\\
[0.5ex]\hline 
& & & & & \\ [-1.5ex]
hyperbolic \& $q=1$ &  & $C^\infty$ \& definable & & $C^\infty$ \& definable & \cite{RainerOmin}\\
[0.5ex]\hline 
& & & & & \\ [-1.5ex]
hyperbolic \& $q=1$  &  & $C^n$ (resp.\ $C^{2n}$) & & $C^1$ (resp.\ twice differentiable) & 
\cite{Bronshtein79}, \cite{Wakabayashi86}, \cite{Mandai85}, \cite{KLM04}, \& \cite{ColombiniOrruPernazza08}\\
[0.5ex]\hline \hline
& & & & & \\ [-1.5ex] 
hyperbolic & & continuous & & continuous (e.g.\ by ordering them increasingly) &  e.g.\ \cite[4.1]{AKLM98}\\
[0.5ex]\hline
& & & & & \\ [-1.5ex]
hyperbolic & &  $C^n$ & & locally Lipschitz & \cite{Bronshtein79} \& \cite{Wakabayashi86} (see also \cite{KurdykaPaunescu08}) \\
[0.5ex]\hline
& & & & & \\ [-1.5ex]
hyperbolic & & $C^\om$ (resp.\ $\cC$, resp.\ arc- & & local desingularization by finitely many & 
\cite{KurdykaPaunescu08} (resp.\ theorem \ref{hypcCperturbth},\\
& & analytic \& subanalytic) & & local blow-ups with smooth center & resp.\ remark \ref{subanrk}) \\
[0.5ex]\hline
& & & & & \\ [-1.5ex]
hyperbolic \& $\tilde \De_s(P)$ has & & $C^\om$ & & locally $C^\om$ & \cite[5.4]{KurdykaPaunescu08}\\
only normal crossings & & & & & \\
[0.5ex]\hline
\end{longtable}
\end{landscape}

}

\section{Preliminaries on polynomials}

\subsection{Coefficients and roots} \label{sturm}
Let
\[
P(z) = z^n + \sum_{j=1}^n (-1)^j a_j z^{n-j} = \prod_{j=1}^n (z-\la_j)
\]
be a univariate monic complex polynomial with coefficients $a_1,\ldots,a_n \in \mathbb{C}$
and roots $\la_1,\ldots,\la_n \in \mathbb{C}$.
By Vieta's formulas, $a_i = \si_i(\la_1,\ldots,\la_n)$, where
$\si_1,\ldots,\si_n$
denote the elementary symmetric functions in $n$ variables:
\begin{equation} \label{esf}
\si_i(\la_1,\ldots,\la_n) =
\sum_{1 \le j_1 < \cdots < j_i \le n} \la_{j_1} \cdots \la_{j_i}.
\end{equation}
It is well-known that each symmetric polynomial in $n$ variables can be written as a polynomial in $\si_1,\ldots,\si_n$, i.e., 
$\C[\la_1,\ldots,\la_n]^{\on{S}_n} = \C[\si_1,\ldots,\si_n]$, where $\on{S}_n$ denotes the symmetric group on $\{1,2,\ldots,n\}$.

Denote by $s_i$ (for $i \in \N$) the Newton polynomials
\begin{equation} \label{Newton}
s_i(\la_1,\ldots,\la_n) = \sum_{j=1}^n \la_j^i
\end{equation}
which are related to the elementary symmetric
functions
by
\begin{equation} \label{rec}
s_k - s_{k-1} \si_1 + s_{k-2} \si_2 - \cdots
+ (-1)^{k-1} s_1 \si_{k-1} + (-1)^k k \si_k = 0, \quad (k \ge 1).
\end{equation}
These relations define a polynomial diffeomorphism $\Ps^n$ such that: 
%We consider the corresponding mappings which are interrelated by a polynomial diffeomorphism $\Ps^n$, given by \eqref{rec}:
\begin{align*}
\si^n &= (\si_1,\ldots,\si_n) : \C^n \to \C^n,\\
s^n &=(s_1,\ldots,s_n) : \C^n \to \C^n,\\
s^n &= \Ps^n \o \si^n. 
\end{align*}
It is easy to compute the Jacobian determinants 
$\det (d s^n(\la)) = n! \prod_{i<j} (\la_j-\la_i)$, 
$\det (d \Ps^n(\si^n))  = (-1)^{n(n-1)/2} n!$,
and, hence, 
\begin{equation} \label{detsi}
\det (d \si^n(\la)) = \prod_{i<j} (\la_i-\la_j).
\end{equation}

Let us consider the so-called Bezoutiant
\[
B :=
\begin{pmatrix}
s_0 & s_1 & \ldots & s_{n-1}\\
s_1 & s_2 & \ldots & s_n \\
\vdots & \vdots & \ddots & \vdots\\
s_{n-1} & s_n & \ldots &  s_{2n-2}
\end{pmatrix}
= \left(s_{i+j-2}\right)_{1 \le i,j \le n}.
\]
Since the entries of $B$ are symmetric polynomials in $\la_1,\ldots,\la_n$,
there exists a unique symmetric $n \times n$ matrix $\tilde{B}$ with
$B = \tilde{B} \circ \si^n$.

Let $B_k$ denote the minor formed by the first $k$ rows and columns of $B$.
Then it is easy to see that
\begin{equation} \label{eqdel}
\Delta_k(\la) := \det B_k(\la) = \sum_{i_1 < i_2 < \cdots < i_k}
(\la_{i_1}-\la_{i_2})^2 \cdots (\la_{i_1}-\la_{i_k})^2 \cdots
(\la_{i_{k-1}}-\la_{i_k})^2.
\end{equation}
In particular, $\De_1(\la)=s_0=n$.
Since the polynomials $\Delta_k$ are symmetric, we have
$\Delta_k = \tilde{\Delta}_k \circ \si^n$
for unique polynomials $\tilde{\Delta}_k$.
By \eqref{eqdel}, the number of
distinct roots of $P$ equals the
maximal $k$ such that $\tilde \De_k(P) \ne 0$.
(Abusing notation we identify $P$ with the $n$-tuple $(a_1,\ldots,a_n)$ of its coefficients when convenient.)

\begin{theorem}[Sylvester's version of Sturm's theorem, e.g.\ {\cite{Procesi78}}] \label{hypth}
Suppose that all coefficients of $P$ are real.
Then all roots of $P$ are real if and only if the symmetric $n\times n$
matrix $\tilde{B}(P)$ is
positive semidefinite. The rank of $\tilde{B}(P)$ equals the number of
distinct roots
of $P$ and its signature equals the number of distinct real roots.
\end{theorem}

\subsection{Hyperbolic polynomials} \label{hyp}
If all roots $\la_j$ (and thus all coefficients $a_j$) of $P$ are real, we say that $P$ is \emph{hyperbolic}.

The space of all hyperbolic polynomials $P$ of fixed degree $n$ can be identified with the semialgebraic subset 
$\si^n(\R^n) \subseteq \R^n$. Its structure is described in theorem \ref{hypth}.
If the roots are ordered increasingly, i.e.,
\[
\la_1(P) \le \la_2(P) \le \cdots \le \la_n(P), \quad \text{ for all } P \in \si^n(\R^n),
\]
then each root $\la_i : \si^n(\R^n) \to \R$ (for $1 \le i \le n$) is continuous (e.g.\ \cite[4.1]{AKLM98}).

Note that all roots of a hyperbolic polynomial $P$ with $a_1=a_2=0$ are
equal to $0$, since
\[
\sum \la_i^2 = s_2(\la)=\si_1(\la)^2-2 \si_2(\la) = a_1^2-2a_2.
\]

Replacing the variable $z$ by $z-a_1(P)/n$ transforms any polynomial $P$ to
another polynomial $\bar P$ with $a_1(\bar P)=0$.
If all roots of $\bar P$ coincide, they have to be equal to $0$. We use that fact repeatedly.

\begin{proposition}[Bounded roots] \label{bdroot}
Let $(P_m)$ be a sequence of univariate monic polynomials over $\C$ with fixed degree $n$ and bounded coefficients. 
If $(\la_m) \subseteq \C$ such that $P_m(\la_m)=0$ for all $m$, then $(\la_m)$ is bounded.
\end{proposition}

\begin{demo}{Proof}
If $a_{m,j}$ denote the coefficients of $P_m$, we find
\begin{equation} \label{bdeq}
|\la_m|^n  \le \sum_{j=1}^n |a_{m,j}||\la_m|^{n-j}.
\end{equation}
Suppose that $(\la_m)$ is unbounded.
Without loss we may assume that $0 < |\la_m| \nearrow \infty$. 
Dividing \eqref{bdeq} by $|\la_m|^{n-1}$ yields a contradiction.
\qed\end{demo}

\section{\texorpdfstring{$C^\infty$}{} classes that admit resolution of singularities} \label{secC}

Following \cite[Section 3]{BM04} we discuss classes of smooth functions that admit resolution of singularities.

\subsection{Classes \texorpdfstring{$\cC$}{} of \texorpdfstring{$C^\infty$}{}-functions} \label{cC}
Let us assume that for every open $U \subseteq \R^q$, $q \in \N$, we have a subalgebra $\cC(U)$ 
of $C^\infty(U)=C^\infty(U,\R)$.
Resolution of singularities in $\cC$ (see \ref{res}) requires only the following assumptions $(\ref{cC}.1)-(\ref{cC}.6)$ on $\cC(U)$, 
for any open $U \subseteq \R^q$.
\begin{enumerate}
\item[(\ref{cC}.1)] $\mathcal P(U) \subseteq \cC(U)$, where $\mathcal P(U)$ denotes the algebra of restrictions to $U$ of polynomial functions on $\R^q$.
\item[(\ref{cC}.2)] {\it $\cC$ is closed under composition}. If $V \subseteq \R^p$ is open and $\vh=(\vh_1,\ldots,\vh_p) : U \to V$ is a mapping 
with each $\vh_i \in \cC(U)$, then $f \o \vh \in \cC(U)$, for all $f \in \cC(V)$.
\end{enumerate}
A mapping $\vh : U \to V$ is called a \emph{$\cC$-mapping} if $f \o \vh \in \cC(U)$, for every $f \in \cC(V)$. 
It follows from (\ref{cC}.1) and (\ref{cC}.2) that $\vh=(\vh_1,\ldots,\vh_p)$ is a $\cC$-mapping if and only if $\vh_i \in \cC(U)$, 
for all $1 \le i \le p$.
\begin{enumerate}
\item[(\ref{cC}.3)] {\it $\cC$ is closed under derivation}. If $f \in \cC(U)$ and $1 \le i \le q$,
then 
$\p_i f \in \cC(U)$. 
\item[(\ref{cC}.4)] {\it $\cC$ is quasianalytic}. If $f \in \cC(U)$ and $\hat f_a=0$, for $a \in U$, then $f$ vanishes in a neighborhood of $a$. 
\end{enumerate}
Since $\{x : \hat f_x=0\}$ is closed in $U$, (\ref{cC}.4) is equivalent to the following property: 
If $U$ is connected, then, for each $a \in U$,
the Taylor series homomorphism $\cC(U) \to \cF_q$, 
$f \mapsto \hat f_a$, is injective.
\begin{enumerate}
\item[(\ref{cC}.5)] {\it $\cC$ is closed under division by a coordinate}. If $f \in \cC(U)$ is identically $0$ along a hyperplane $\{x : x_i=a_i\}$, i.e.,
$f(x_1,\ldots,x_{i-1},a_i,x_{i+1},\ldots,x_q) \equiv 0$, then $f(x)= (x_i-a_i) h(x)$, where $h \in \cC(U)$.
\item[(\ref{cC}.6)] {\it $\cC$ is closed under taking the inverse}. 
Let $\vh : U \to V$ be a $\cC$-mapping between open subsets $U$ and $V$ in $\R^q$. 
Let $a \in U$, $\vh(a)=b$, and suppose that the Jacobian matrix $(\p \vh/\p x)(a)$ is invertible. Then 
there exist neighborhoods $U'$ of $a$, $V'$ of $b$, and a $\cC$-mapping $\ps : V' \to U'$ such that $\ps(b)=a$ and $\vh \o \ps = \on{id}_{V'}$. 
\end{enumerate}
Property (\ref{cC}.6) is equivalent to the \emph{implicit function theorem in $\cC$}:
Let $U \subseteq \R^q \times \R^p$ be open. Suppose that $f_1,\ldots,f_p \in \cC(U)$, $(a,b) \in U$, $f(a,b)=0$, 
and $(\p f/\p y)(a,b)$ is invertible, where $f=(f_1,\ldots,f_p)$.
Then there is a neighborhood $V \times W$ of $(a,b)$ in $U$ and a $\cC$-mapping $g : V \to W$ such that $g(a)=b$ and $f(x,g(x))=0$, for $x \in V$.

It follows from (\ref{cC}.6) that {\it $\cC$ is closed under taking the reciprocal}: If $f \in \cC(U)$ vanishes nowhere in $U$, 
then $1/f \in \cC(U)$.

A complex valued function $f : U \to \C$ is said to be a \emph{$\cC$-function}, or to belong to $\cC(U,\C)$, 
if $(\Re f, \Im f) : U \to \R^2$ is a $\cC$-mapping.
It is immediately verified that (\ref{cC}.3)--(\ref{cC}.5) hold for complex valued functions $f \in \cC(U,\C)$ as well.

\textbf{From now on, unless otherwise stated, 
$\cC$ will denote a fixed, but arbitrary, class of $C^\infty$-functions satisfying the conditions (\ref{cC}.1)-(\ref{cC}.6).}

\begin{lemma}[Splitting lemma in $\cC$, cf.\ {\cite[3.4]{AKLM98}}] \label{split}
Let $P_0 = z^n + \sum_{j=1}^n (-1)^j a_j z^{n-j}$ be a complex polynomial
satisfying
$P_0 = P_1 \cdot P_2$, where $P_1$ and $P_2$ are monic polynomials without common
root.
Then for $P$ near $P_0$ we have $P = P_1(P) \cdot P_2(P)$
for $\cC$-mappings
of monic polynomials $P \mapsto P_1(P)$ and $P \mapsto P_2(P)$,
defined for $P$
near $P_0$, with the given initial values.
(Here $P \mapsto P_i(P)$ is understood as a mapping $\R^{2n} \to \R^{2 \deg P_i}$.)
\end{lemma}

\begin{demo}{Proof}
Let the polynomial $P_0$ be represented as the product
\[
P_0 = P_1 \cdot P_2 
= \Big(z^p + \sum_{j=1}^p (-1)^j b_j z^{p-j}\Big) \cdot \Big(z^q + \sum_{j=1}^q (-1)^j c_j z^{q-j}\Big),
\]
where $p + q = n$. Let $\la_1,\ldots,\la_n$ be the roots of $P_0$, ordered in such a way 
that the first $p$ are the roots of $P_1$ and the last $q$ are those of $P_2$. 
There is a polynomial mapping $\Ph^{p,q}$ such that
$(a_1,\ldots,a_n) = \Ph^{p,q}(b_1,\ldots,b_p,c_1,\ldots,c_q)$. 
Let $b = (b_1,\ldots,b_p)$ and $c = (c_1,\ldots,c_q)$.
Then
\begin{gather*}
\si^n = \Phi^{p,q} \circ (\si^p \times \si^q),\\
\det(d \si^n) = \det(d \Ph^{p,q} (b,c)) \det(d \si^p) \det(d \si^q),
\end{gather*}
and, by \eqref{detsi}, 
\[
\det(d \Ph^{p,q}(b,c)) = \prod_{1 \le i \le p < j \le n} (\la_i - \la_j) \ne 0,
\]
since $P_1$ and $P_2$ do not have common roots. 

If we view $\Ph^{p,q}$ as a mapping $\R^{2n} \to \R^{2n}$, then its Jacobian determinant at $(b,c)$ is still not $0$,
by lemma \ref{RCdet} below.
So, by (\ref{cC}.1) and (\ref{cC}.6), 
$\Ph^{p,q}$ is a $\cC$-diffeomorphism near $(b,c)$.
\qed\end{demo}

\begin{lemma} \label{RCdet}
Let $A=(A_{ij}) \in \C^{n \times n}$. Consider the block matrix $B=(B_{ij}) \in \R^{2n \times 2n}$, where
\[
B_{ij} = \begin{pmatrix}
\on{Re} A_{ij} & - \on{Im} A_{ij}\\ \on{Im} A_{ij} & \on{Re} A_{ij}
\end{pmatrix}, 
\quad (1 \le i,j \le n).
\]
Then $\det_{\R} B = |\det_{\C} A|^2$.
\qed\end{lemma}

\subsection{\texorpdfstring{$\cC$}{}-manifolds} \label{Cmf}

One can use the open subsets $U \subseteq \R^q$ and the 
algebras of functions $\cC(U)$ as local models to define a category $\underline{\cC}$ of \emph{$\cC$-manifolds} and 
\emph{$\cC$-mappings}. The dimension theory of $\underline{\cC}$ follows from that of $C^\infty$-manifolds.

The implicit function property (\ref{cC}.6) implies that a \emph{smooth} (not singular) subset of a $\cC$-manifold is a $\cC$-submanifold:

\begin{proposition}
Let $M$ be a $\cC$-manifold. Suppose that $U$ is open in $M$, $g_1,\ldots,g_p \in \cC(U)$, and the gradients 
$\nabla g_i$ are linearly independent at every point of the zero set $X:=\{x \in U : g_i(x)=0 \text{ for all }i\}$. 
Then $X$ is a closed $\cC$-submanifold of $U$ of codimension $p$. \qed
\end{proposition}

\section{Quasianalytic Denjoy--Carleman classes} \label{secDC}

\subsection{Denjoy--Carleman classes} \label{DC} 
See \cite{Thilliez08} and references therein.
Let $U \subseteq \R^q$ be open. 
Let $M=(M_k)_{k \in \N}$ be a non-decreasing sequence of real numbers with $M_0=1$.
We denote by $C^M(U)$ the set of all $f \in C^\infty(U)$ such that for every compact $K \subseteq U$ there are constants $C,\rh >0$ with 
\begin{equation} \label{CM}
|\p^\al f(x)| \le C \rh^{|\al|} |\al|!\, M_{|\al|} \quad \text{ for all } \al \in \N^q \text{ and } x \in K.
\end{equation} 
We call $C^M(U)$ a \emph{Denjoy--Carleman class} of functions on $U$.
If $M_k=1$, for all $k$, then $C^M(U)$ coincides with the ring $C^\om(U)$ of real analytic functions
on $U$. In general, $C^\om(U) \subseteq C^M(U) \subseteq C^\infty(U)$. 
Hence $\cC=C^M$ satisfies property (\ref{cC}.1).

We assume that $M=(M_k)$ is \emph{logarithmically convex}, i.e.,
\begin{equation} \label{logconvex}
M_k^2 \le M_{k-1} \, M_{k+1} \quad \text{ for all } k,
\end{equation}
or, equivalently, $M_{k+1}/M_k$ is increasing.
Using $M_0=1$, we obtain that also $(M_k)^{1/k}$ is increasing and
\begin{equation}
M_l \, M_k\le M_{l+k} \quad \text{  for all }l,k\in \N.
\label{logconvex1}
\end{equation}

Hypothesis \eqref{logconvex} implies that $C^M(U)$ is a ring, for all open
subsets $U \subseteq \R^q$, 
which can easily be derived from \eqref{logconvex1} by means of Leibniz' rule.
Note that definition \eqref{CM} makes sense also for mappings
$U\to \mathbb R^p$.
For $C^M$-mappings, \eqref{logconvex} guarantees stability under composition
(\cite{Roumieu62/63}, \cite[4.7]{BM04}). So $\cC=C^M$ satisfies property (\ref{cC}.2).

A further consequence of \eqref{logconvex} is the inverse function theorem
for $C^M$ (\cite{Komatsu79}, \cite[4.10]{BM04}). Thus $\cC=C^M$ satisfies property (\ref{cC}.6).

Suppose that $M=(M_k)$ and $N=(N_k)$ satisfy %$M_k \le C^k \, N_k$, for all $k$ and a constant $C$, or equivalently,
\begin{equation} \label{incl}
\sup_{k \in \N_{>0}} \Big(\frac{M_k}{N_k}\Big)^{\frac{1}{k}} < \infty.
\end{equation}
Then, evidently $C^M(U) \subseteq C^N(U)$. The converse
is true as well: %(if \eqref{logconvex} is assumed; actually the weaker assumption that $k!M_k$ is logarithmically convex suffices): 
There exists $f \in C^M(\R)$ such that 
$|f^{(k)}(0)| \ge k! \, M_k$ for all $k$ (see \cite[Theorem 1]{Thilliez08}).
So the inclusion $C^M(U) \subseteq C^N(U)$ implies \eqref{incl}. 

Setting $N_k=1$ in \eqref{incl} yields that $C^\om(U) = C^M(U)$ if and only if 
\[
\sup_{k \in \N_{>0}} (M_k)^{\frac{1}{k}} < \infty.
\] 
As $(M_k)^{1/k}$ is increasing (by \eqref{logconvex}), 
the strict inclusion $C^\om(U) \subsetneq C^M(U)$ is equivalent to  
\[
\lim_{k \to \infty} (M_k)^{\frac{1}{k}} = \infty.
\] 

The class $\cC=C^M$ is stable under derivation (property (\ref{cC}.3)) if and only if
\begin{equation} \label{der}
\sup_{k \in \N_{>0}} \Big(\frac{M_{k+1}}{M_k}\Big)^{\frac{1}{k}} < \infty.
\end{equation}
The first order partial derivatives of elements in $C^M(U)$
belong to $C^{M^{+1}}(U)$,
where $M^{+1}$ denotes the shifted sequence 
$M^{+1} = (M_{k+1})_{k \in \N}$.
So the equivalence follows from \eqref{incl}, by replacing $M$ with
$M^{+1}$ and $N$ with $M$.

By the standard integral formula, stability under derivation implies that $\cC=C^M$ fulfills property (\ref{cC}.5).

\subsection{Quasianalyticity}
Suppose that $M$ is logarithmically convex (actually, logarithmic convexity of $k!M_k$ suffices). Then,
by the Denjoy--Carleman theorem (\cite{Denjoy21}, \cite{Carleman26}), 
$\cC=C^M$ is quasianalytic (satisfies (\ref{cC}.4)) if and only if 
\begin{equation} \label{eqqa}
\sum_{k=1}^\infty \frac1{(k!M_k)^{1/k}} = \infty \quad \text{or, equivalently, } \quad \sum_{k=0}^\infty \frac{M_k}{(k+1)M_{k+1}}=\infty. 
\end{equation}
For contemporary proofs see for instance \cite[1.3.8]{Hoermander83I} 
or \cite[19.11]{Rudin87}.

\begin{proposition}
If $M$ is a non-decreasing sequence of real numbers with $M_0=1$
satisfying \eqref{logconvex}, \eqref{der}, and \eqref{eqqa},
then the Denjoy--Carleman class $\cC=C^M$ has the properties \emph{(\ref{cC}.1)--(\ref{cC}.6)}. 
If $C^M$ is not closed under derivation (i.e., \eqref{der} fails), 
then $\cC = \bigcup_{j \in \N} C^{M^{+j}}$ has the properties \emph{(\ref{cC}.1)--(\ref{cC}.6)}.\qed
\end{proposition}

\section{Resolution of singularities in \texorpdfstring{$\cC$}{}}

\subsection{Blow-ups} \label{blowup}
Let $M$ be a smooth manifold and let $C$ be a smooth closed subset of $M$. 
The \emph{blow-up} of $M$ with \emph{center} $C$ is a proper smooth mapping $\vh : M' \to M$
from a smooth manifold $M'$ that can be described in local coordinates as follows. 

Let $U \subseteq \R^q$ be an open neighborhood of $0$ and let $C = \{x_i=0 ~\text{for}~ i \in I\}$ be a 
coordinate subspace, where $I$ is a subset of $\{1,\ldots,q\}$. 
The \emph{blow-up} $\vh : U' \to U$ with \emph{center} $C$ is a mapping 
where $U'$ can be covered by coordinate charts $U_i'$, for $i \in I$, and each $U_i'$ has a coordinate system $y_1,\ldots,y_q$ in which $\vh$ is given by 
\begin{align*}
x_j = \left\{
\begin{array}{l@{\quad \text{for} ~}l}
y_i, & j=i \\
y_i y_j, & j \in I \setminus\{i\} \\
y_j, & j \not\in I
\end{array}
\right..
\end{align*}
Assuming (without loss) $I=\{1,\ldots,p\}$ and $x = (\bar x,\tilde x) \in \R^p \times \R^{q-p}$, we have
\[
U' \cong \{(x,\xi) \in U \times \mathbb{RP}^{p-1} : \bar x \in \xi\},
\]
and, if we use homogeneous coordinates $\xi = [\xi_1,\ldots,\xi_p]$,
\[
U'= \{(x,\xi) \in U \times \mathbb{RP}^{p-1} : x_i \xi_j = x_j \xi_i \text{ for } 1 \le i,j \le p\}.
\]
We can cover $U'$ by coordinate charts $U_i'= \{(x,\xi) \in U' : \xi_i \ne 0\}$, for $i \in I$, with coordinates $y_1,\ldots,y_q$ where
\begin{align*}
y_j = \left\{
\begin{array}{l@{\quad \text{for} ~}l}
x_i, & j=i \\
\frac{\xi_j}{\xi_i}, & j \in I \setminus\{i\} \\
x_j, & j \not\in I
\end{array}
\right..
\end{align*}

The \emph{blow-up} of a smooth manifold $M$ with \emph{center} a smooth closed subset $C$ is a smooth mapping 
$\vh : M' \to M$, where $M'$ is a smooth manifold, such that: 
\begin{enumerate}
\item Every point of $C$ admits a coordinate neighborhood $U$ in which $C$ is a coordinate subspace and over $U$ the 
mapping $\vh : M' \to M$ identifies with the mapping $U' \to U$ from above.
\item $\vh$ restricts to a diffeomorphism over $M \setminus C$.
\end{enumerate}
These conditions determine $\vh : M' \to M$ uniquely up to a diffeomorphism of $M'$ commuting with $\vh$.
If $\on{codim} C = 1$ then the blow-up $\vh$ is the identity.

If $M$ is a $\cC$-manifold and $\vh : M' \to M$ is the blow-up with center a closed $\cC$-submanifold $C$ of $M$, 
then $M'$ is a $\cC$-manifold and $\vh$ is a $\cC$-mapping (cf.\ \cite[3.9]{BM04}):

\begin{proposition}
The category $\underline \cC$ of $\cC$-manifolds and $\cC$-mappings is closed under blowing up with center a closed $\cC$-submanifold. \qed 
\end{proposition}

\subsection{Resolution of singularities} \label{res}

We shall use a simple version of the desingularization theorem of Hironaka \cite{Hironaka64} for $\cC$-function classes due 
to Bierstone and Milman \cite{BM04}.
We use the terminology therein.

Let us regard a $\cC$-manifold $M$ as local-ringed space $(|M|,\cO^\cC_M)$ with $|M|$ the underlying 
topological space of $M$ and $\cO^\cC_M$ the sheaf of germs of $\cC$-functions at points of $M$. 
Let $\cI \subseteq \cO^\cC_M$ be a sheaf of ideals of finite type, i.e., for each $a \in M$, there is an open neighborhood $U$ of $a$ and finitely many sections 
$f_1,\ldots,f_p \in \cO^\cC_M(U)=\cC(U)$ such that, for all $b \in U$, the stalk $\cI_b$ is generated by the germs of the $f_i$ at $b$.
Put $|X|:=\on{supp} \cO^\cC_M/\cI$ and 
$\cO^\cC_X := (\cO^\cC_M/\cI)|_{|X|}$. Then $X=(|X|,\cO^\cC_X)$ is called a \emph{closed $\cC$-subspace} of $M$, 
and we write $\cI=\cI_X$. It is a \emph{hypersurface} if $\cI_X$ is a sheaf of 
principal ideals. A closed $\cC$-subspace $X$ is \emph{smooth} at $a \in X$ if $\cI_{X,a}$ is generated by elements 
with linearly independent gradients at $a$. By proposition \ref{Cmf}, a smooth $\cC$-subspace is a $\cC$-submanifold.

Let $\vh : N \to M$ be a $\cC$-mapping of $\cC$-manifolds. If $\cI \subseteq \cO^\cC_M$ is a sheaf of ideals of finite type, we denote by $\vh^{-1}(\cI) \subseteq \cO^\cC_N$ the ideal sheaf $\vh^*(\cI) \cdot \cO^\cC_N$ 
whose stalk at each $b \in N$ is generated by the ring of pullbacks $\vh^*(\cI)_b$ of all elements in $\cI_{\vh(b)}$.
If $X$ is a closed $\cC$-subspace of $M$, let $\vh^{-1}(X)$ denote the closed $\cC$-subspace of $N$ determined by 
$\vh^{-1}(\cI_X)$.

Let $M$ be a $\cC$-manifold, $C$ a $\cC$-submanifold of $M$, and let $\vh : M' \to M$ be the blow-up of $M$ 
with center $C$. Then $\vh^{-1}(C)$ is a smooth closed subspace in $M'$. We denote by $y_{\on{exc}}$ a 
generator of $\cI_{\vh^{-1}(C),a'}$, at any $a' \in M'$. 

Let $X \subseteq M$ be a hypersurface. 
The \emph{strict transform $X'$ of $X$ by $\vh$} is the hypersurface of $M'$ determined by $\cI_{X'}$, where 
$\cI_{X'} \subseteq \cO^\cC_{M'}$ is defined as follows: If $a' \in M'$, $a = \vh(a')$, and $g$ is a generator of 
$\cI_{X,a}$, then $\cI_{X',a'}$ is the ideal generated by $g':= y_{\on{exc}}^{-d} g \o \vh$, where $d$ is the 
largest power of $y_{\on{exc}}$ that factors from $g \o\vh$. 
(If $a' \not\in \vh^{-1}(C)$, then we may take $y_{\on{exc}}=1$.)  
See \cite[5.6]{BM04} and \cite[Section 3]{BM97} for the difference between \emph{weak} and \emph{strict} transform 
(and the problems with the latter in $\cC$) if $X$ is \emph{not} a hypersurface. 

We say that a hypersurface $X$ has only \emph{normal crossings}, if locally there exist suitable coordinates in which 
$\cI_X$ is generated by a monomial.

\begin{theorem}[{\cite[5.12]{BM04}}]
\label{resth}
Let $M$ be a $\cC$-manifold, $X$ a closed $\cC$-hypersurface in $M$, and $K$ a compact subset of $M$. 
Then, there is a neighborhood $W$ of $K$ and a surjective mapping $\vh : W' \to W$ of class $\cC$, such that:
\begin{enumerate}
\item $\vh$ is a composite of finitely many $\cC$-mappings, 
each of which is either a blow-up with smooth center 
(that is nowhere dense in the smooth points of the strict transform of $X$) 
or a surjection of the form $\bigsqcup_j U_j \to \bigcup_j U_j$, 
where the latter is a finite covering of the target space by coordinate charts.
\item The final strict transform $X'$ of $X$ is smooth, and $\vh^{-1}(X)$ has only normal crossings. 
(In fact $\vh^{-1}(X)$ and $\det d \vh$ simultaneously have only normal crossings, 
where $d \vh$ is the Jacobian matrix of $\vh$ with respect to any local coordinate system.)
\end{enumerate}
\end{theorem}

See \cite[5.9 \& 5.10]{BM04} and \cite{BM97} for stronger desingularization theorems in $\cC$. 

\section{Quasianalytic perturbation of polynomials}

We prove in this section that the roots of a $\cC$-family of polynomials $P$ can be parameterized locally by $\cC$-functions 
after modifying $P$ in a precise way.

\subsection{Local blow-ups and local power substitutions} \label{not}
We introduce notation following \cite[Section 4]{BM88}.

Let $M$ be a $\cC$-manifold. 
A family of $\cC$-mappings $\{\pi_j : U_j \to M\}$ is called a \emph{locally finite covering} of $M$ if
the images $\pi_j(U_j)$ are subordinate to a locally finite open covering $\{W_j\}$ of $M$ (i.e.\ $\pi_j(U_j) \subseteq W_j$ for all $j$) and if, 
for each compact $K \subseteq M$, there are compact $K_j \subseteq U_j$ such that $K = \bigcup_j \pi_j(K_j)$ (the union is finite). 

Locally finite coverings can be \emph{composed} in the following way (see \cite[4.5]{BM88}):
Let $\{\pi_j : U_j \to M\}$ be a locally finite covering of $M$, and let $\{W_j\}$ be as above.
For each $j$, suppose that $\{\pi_{ji} : U_{ji} \to U_j\}$ is a 
locally finite covering of $U_j$. We may assume without loss that the $W_j$ are relatively compact.
(Otherwise, choose a locally finite covering $\{V_j\}$ of $M$ by relatively compact open subsets. 
Then the mappings $\pi_j|_{\pi_j^{-1}(V_i)} :  \pi_j^{-1}(V_i) \to M$, for all $i$ and $j$, 
form a locally finite covering of $M$.)
Then, for each $j$, there is a finite subset $I(j)$ of $\{i\}$ such that the $\cC$-mappings 
$\pi_j \o \pi_{ji} : U_{ji} \to M$, for all $j$ and all $i \in I(j)$, form a locally finite covering of $M$. 

We shall say that $\{\pi_j\}$ is a \emph{finite covering}, if $j$ varies in a finite index set.

A \emph{local blow-up $\Ph$} over an open subset
$U$ of $M$ means the composition $\Ph = \iota \o \vh$ of a blow-up $\vh : U' \to U$ with smooth center and 
of the inclusion $\iota : U \to M$. 

We denote by \emph{local power substitution} a mapping of $\cC$-manifolds $\Ps: V \to M$ of the form 
$\Ps = \iota \o \ps$, where $\iota : W \to M$ is the inclusion of a coordinate chart $W$ of $M$ and 
$\ps : V \to W$ is given by 
\begin{equation} \label{defps}
(y_1,\ldots,y_q) = \ps_{\ga,\ep}(x_1,\ldots,x_q) := ((-1)^{\ep_1} x_1^{\ga_1},\ldots,(-1)^{\ep_q} x_q^{\ga_q}),
\end{equation}
for some $\ga=(\ga_1,\ldots,\ga_q) \in (\N_{>0})^q$ and all $\ep = (\ep_1,\ldots,\ep_q) \in \{0,1\}^q$, 
where $y_1,\ldots,y_q$ denote the coordinates of $W$ (and $q = \dim M$).  

\subsection{}

We consider the natural partial ordering of multi-indices: 
If $\al,\be \in \N^q$, then $\al \le \be$ means $\al_i \le \be_i$ for all $1 \le i \le q$.

\begin{lemma}[{\cite[7.7]{BM04}} or {\cite[4.7]{BM88}}] \label{order} 
Let $\al,\be,\ga \in \N^q$ and let $a(x),b(x),c(x)$ be non-vanishing germs of real or complex valued functions of class $\cC$ at the origin of 
$\R^q$. If 
\[
x^\al a(x) - x^\be b(x) = x^\ga c(x),
\]
then either $\al \le \be$ or $\be \le \al$.
\end{lemma}

\begin{demo}{Proof}
Let $\de = (\de_1,\dots,\de_q)$ where $\de_k=\min\{\al_k,\be_k\}$. If $\de=\al$ then $\al \le \be$. Otherwise, $\de_k \ne \al_k$ for some $k$. 
On $\{x_k=0\}$ we have $x^{\al-\de} = 0$ and $0 \ne - x^{\be-\de} b(x) = x^{\ga-\de} c(x)$. 
Since $b$ and $c$ are non-vanishing, we obtain $\be=\ga$, by (\ref{cC}.5). %and (\ref{cC}.6). 
So $x^\al a(x) = x^\be (b(x)+c(x))$ and hence $\al \ge \be$, again by (\ref{cC}.5).
\qed\end{demo}

\subsection{} \label{ncross}
Let $M$ be a $\cC$-manifold and let $f$ be a real or complex valued $\cC$-function on $M$. 
We say that $f$ has only \emph{normal crossings} if each point in $M$ admits a coordinate neighborhood $U$ with coordinates 
$x=(x_1,\ldots,x_q)$ such that
\[
f(x)=x^\al g(x), \quad x \in U,
\]
where $g$ is a non-vanishing $\cC$-function on $U$, and $\al \in \N^q$.
Observe that, if a product of functions has only normal crossings, then each factor has only normal crossings.
For: Let $f_1, f_2,g$ be $\cC$-functions defined near $0 \in \R^q$ such that $f_1(x) f_2(x)=x^\al g(x)$ and $g$ is non-vanishing.
By quasianalyticity (\ref{cC}.4), $f_1 f_2|_{\{x_j=0\}}=0$ implies $f_1|_{\{x_j=0\}}=0$ or $f_2|_{\{x_j=0\}}=0$.
So the assertion follows from (\ref{cC}.5).

\subsection{} \label{RtoC}
Let $M$ be a $\cC$-manifold, $K \subseteq M$ be compact, and $f \in \cC(M,\C)$. Then the exists a neighborhood $W$ of $K$ 
and a finite covering $\{\pi_k : U_k \to W\}$ of $W$ by $\cC$-mappings $\pi_k$,
each of which is a composite of finitely many local blow-ups with smooth center, such that, for each $k$, 
the function $f \o \pi_k$ has only normal crossings. This follows from theorem \ref{resth} applied to the real valued $\cC$-function 
$|f|^2=f\overline f$ and the observation in \ref{ncross}.

\subsection{Reduction to smaller permutation groups} \label{reduce}
In the proof of theorem \ref{cCperturb} we shall reduce our perturbation
problem in virtue of the splitting lemma \ref{split}:

The space $\on{Pol}^n$ of polynomials $P(z)=z^n+\sum_{j=1}^n (-1)^j a_j
z^{n-j}$ of fixed degree $n$
naturally identifies with $\C^n$ (by mapping $P$ to $(a_1,\ldots,a_n)$).
Moreover, $\on{Pol}^n$ may be
viewed as the orbit space $\C^n/\on{S}_n$ with respect to the standard action of the
symmetric group $\on{S}_n$ on
$\C^n$ by permuting the coordinates (the roots of $P$).
In this picture the mapping $\si^n : \C^n \to \C^n$ %(introduced in \ref{sturm}) 
identifies with the orbit
projection $\C^n \to \C^n/\on{S}_n$, since the elementary symmetric
functions $\si_i$ in \eqref{esf} generate
the algebra of symmetric polynomials on $\C^n$, i.e., $\C[\C^n]^{\on{S}_n}=\C[\si_1,\ldots,\si_n]$.

Consider a family of polynomials
\[
P(x)(z) = z^n + \sum_{j=1}^n (-1)^j a_j(x) z^{n-j},
\]
where the coefficients $a_j$ are complex valued $\cC$-functions defined in a $\cC$-manifold $M$. 
Let $x_0 \in M$.
If $P(x_0)$ has distinct roots $\nu_1,\ldots,\nu_l$,  
the splitting lemma \ref{split} provides a $\cC$-factorization
$P(x) = P_1(x) \cdots P_l(x)$ near $x_0$ such that no two factors have common roots and
all roots of $P_h(x_0)$ are
equal to $\nu_h$, for $1 \le h \le l$.
This factorization amounts
to a reduction of
the $\on{S}_n$-action on $\C^n$ to the $\on{S}_{n_1} \times \cdots \times
\on{S}_{n_l}$-action
on $\C^{n_1} \oplus \cdots \oplus \C^{n_l}$, where $n_h$ is the
multiplicity of $\nu_h$.

We shall use the following notation:
\[
\on{S}(P(x_0)) := \on{S}_{n_1} \times \cdots \times \on{S}_{n_l},
\]
iff $P(x_0)$ has $l$ pairwise distinct roots with respective multiplicities $n_1,\dots,n_l$.

Further, we will \emph{remove fixed points} of the 
$\on{S}_{n_1} \times \cdots \times \on{S}_{n_l}$-action on
$\C^{n_1} \oplus \cdots \oplus \C^{n_l}$ or, equivalently, reduce each
factor
$P_h(x)(z) = z^{n_h} + \sum_{j=1}^{n_h} (-1)^j a_{h,j}(x) z^{n_h-j}$ to the
case $a_{h,1}=0$ by replacing $z$
by $z-a_{h,1}(x)/n_h$. The effect on the roots of $P_h$ is a shift by a $\cC$-function.

If $P$ is hyperbolic, we consider the $\on{S}_n$-module $\R^n$ instead of
$\C^n$. In that case the orbit space $\R^n/\on{S}_n$ identifies with the semialgebraic subset $\si^n(\R^n) \subseteq \R^n$, whose 
structure is described in theorem \ref{sturm}. 
Evidently, the splitting lemma \ref{split} produces a $\cC$-factorization 
$P = P_1 \cdots P_l$, where each factor $P_h$ is hyperbolic again.

\begin{theorem}[$\cC$-perturbation of polynomials]  \label{cCperturb}
Let $M$ be a $\cC$-manifold.  
Consider a family of polynomials 
\[
P(x)(z) = z^n + \sum_{j=1}^n (-1)^j a_j(x) z^{n-j}, 
\]
with coefficients $a_j$ (for $1 \le j \le n$) in $\cC(M,\C)$.
Let $K \subseteq M$ be compact.
Then there exist:
\begin{enumerate}
\item a neighborhood $W$ of $K$, and
\item a finite covering $\{\pi_k : U_k \to W\}$ of $W$, where each  
$\pi_k$ is a composite of finitely many mappings each of which is either a 
local blow-up with smooth center or a local power substitution (in the sense of \ref{not}),
\end{enumerate}
such that, for all $k$, the family of polynomials 
$P \o \pi_k$ allows a $\cC$-parameterization of its roots on $U_k$, i.e., there exist
$\la_i^{k} \in \cC(U_k,\C)$ (for $1 \le i \le n$) such that 
\[
P(\pi_k(x))(z) = z^n + \sum_{j=1}^n (-1)^j a_j(\pi_k(x)) z^{n-j} = 
\prod_{i=1}^n (z-\la_i^{k}(x)).
\]
\end{theorem}

\begin{demo}{Proof}
Since the statement is local, we may assume without loss that $M$ is an open neighborhood of $0 \in \R^q$.
In view of \ref{reduce}, we use induction on $|\on{S}(P(0))|$, the order of the permutation group acting on the roots of $P(0)$.

If $|\on{S}(P(0))|=1$, all roots of $P(0)$ are pairwise different.
Then the roots of $P$ may be parameterized in a $\cC$-way near $0$,
by the implicit function theorem (property (\ref{cC}.6)) or by the splitting lemma \ref{split}. 

Suppose that $|\on{S}(P(0))|>1$. Let $\nu_1,\ldots,\nu_l$ denote the distinct roots of $P(0)$; some of them are multiple ($l=1$ is allowed).
The
splitting lemma \ref{split} provides a $\cC$-factorization
$P(x) = P_1(x) \cdots P_l(x)$ near $0$, where, for $1 \le h \le l$,
\[
P_h(x)(z) = z^{n_h} + \sum_{j=1}^{n_h} (-1)^j a_{h,j}(x) z^{n_h-j},
\]
such that no two factors have common roots and all roots of $P_h(0)$ are equal to $\nu_h$. 
As indicated in \ref{reduce}, we reduce to
the $\on{S}_{n_1} \times \cdots \times \on{S}_{n_l}$-action
on $\C^{n_1} \oplus \cdots \oplus \C^{n_l}$ 
and we remove fixed points. So we may assume that $a_{h,1}=0$ for all $h$.

Then all roots of $P_h(0)$ are equal to $0$, and so $a_{h,j}(0) = 0$, for all $1 \le h\le l$ and $1 \le j \le n_h$ (by Vieta's formulas).
If all coefficients $a_{h,j}$ (for $1 \le j \le n_h$) of $P_h$ are identically
$0$, we choose its roots $\la_{h,j}=0$ for all $1 \le j \le n_h$ and remove the factor $P_h$ from the product $P_1 \cdots P_l$.
So we can assume that for each $1 \le h \le l$ there is a $2 \le j \le n_h$ such that $a_{h,j} \ne 0$.

Let us define the $\cC$-functions
\begin{equation} \label{Aa}
A_{h,j}(x) = a_{h,j}(x)^{\frac{n!}{j}}, \quad (\text{for } 1 \le h \le l \text{ and } 2 \le j \le n_h).
\end{equation}
By theorem \ref{resth} (and \ref{RtoC}), we find a finite covering $\{\pi_k : U_k \to U\}$ of a neighborhood $U$ of $0$ 
by $\cC$-mappings $\pi_k$, each of
which is a composite of finitely many local blow-ups with smooth center, such that, for each $k$, 
the non-zero $A_{h,j} \o \pi_k$ (for $1 \le h \le l$ and $2 \le j \le n_h$) and its pairwise non-zero differences 
$A_{h,i} \o \pi_k - A_{m,j} \o \pi_k$ (for $1 \le h \le m \le l$, $1 \le i \le n_h$, and $1 \le j \le n_m$) 
simultaneously have only normal crossings. 

Let $k$ be fixed and let $x_0 \in U_k$. 
Then $x_0$ admits a neighborhood $W_k$ with suitable coordinates in which $x_0=0$ and such that 
(for $1 \le h \le l$ and $2 \le j \le n_h$) 
either $A_{h,j} \o \pi_k=0$ or 
\[
(A_{h,j} \o \pi_k)(x)=x^{\al_{h,j}} A_{h,j}^{k}(x), %\quad \text{ where } A_{h,j}^{k}(0)\ne0.
\]
where $A_{h,j}^{k}$ is a non-vanishing $\cC$-function on $W_k$, and $\al_{h,j} \in \N^q$.
The collection of the multi-indices $\{\al_{h,j} : A_{h,j} \o \pi_k \ne 0, 1 \le h \le l, 2 \le j \le n_h\}$ is totally ordered, by lemma \ref{order}. 
Let $\al$ denote its minimum. 

If $\al=0$, then $(A_{h,j} \o \pi_k)(x_0)=A_{h,j}^k(x_0)\ne 0$ for some $1 \le h \le l$ and $2 \le j \le n_h$. 
So, by \eqref{Aa}, not all roots of $(P_h \o \pi_k)(x_0)$ coincide (since $a_{h,1} \o \pi_k = 0$).
Thus, $|\on{S}((P \o \pi_k)(x_0))| < |\on{S}(P(0))|$, and, by the induction hypothesis, 
there exists a finite covering $\{\pi_{kl} : W_{kl} \to W_k\}$ of $W_k$ (possibly shrinking $W_k$) of the type described in $(2)$ such that, 
for all $l$, the family of polynomials $P \o \pi_k \o \pi_{kl}$ allows a $\cC$-parameterization of its roots on $W_{kl}$.

Let us assume that $\al \ne 0$.
Then there exist $\cC$-functions $\tilde A_{h,j}^{k}$ (some of them $0$) such that, for all $1 \le h \le l$ and $2 \le j \le n_h$,
\begin{equation} \label{Atilde}
(A_{h,j} \o \pi_k)(x)=x^{\al} \tilde A_{h,j}^{k}(x),
\end{equation}
and $\tilde A_{h,j}^{k}(x_0)= A_{h,j}^{k}(x_0)\ne 0$ for some $1 \le h \le l$ and $2 \le j \le n_h$. Let us write
\[
\frac{\al}{n!} = \left(\frac{\al_1}{n!},\ldots,\frac{\al_q}{n!}\right) 
= \left(\frac{\be_1}{\ga_1},\ldots,\frac{\be_q}{\ga_q}\right), 
\]
where $\be_i,\ga_i \in \N$ are relatively prime (and $\ga_i>0$), for all $1 \le i \le q$. 
Put $\be=(\be_1,\ldots,\be_q)$ and $\ga=(\ga_1,\ldots,\ga_q)$. 
Then (by \eqref{Aa} and \eqref{Atilde}), for each $1 \le h \le l$, $2 \le j \le n_h$, and $\ep \in \{0,1\}^q$,
the $\cC$-function $a_{h,j} \o \pi_k \o \ps_{\ga,\ep}$ is divisible by $x^{j \be}$ %, for $1 \le h \le l$ and $2 \le j \le n_h$ 
(where $\ps_{\ga,\ep}$ is defined by \eqref{defps}).
By (\ref{cC}.5), 
there exist $\cC$-functions $a_{h,j}^{k,\ga,\ep}$ such that
\[
(a_{h,j} \o \pi_k \o \ps_{\ga,\ep})(x) = x^{j \be} a_{h,j}^{k,\ga,\ep}(x), \quad (\text{for } 1 \le h \le l \text{ and } 2 \le j \le n_h).
\]
By construction, for some $1 \le h \le l$ and $2 \le j \le n_h$, we have $a_{h,j}^{k,\ga,\ep}(0) \ne 0$, independently of $\ep$.
So there exist a local power substitution $\ps_k : V_k \to W_k$ given in local coordinates by $\ps_{\ga,\ep}$ (for $\ep \in \{0,1\}^q$) 
and functions $a_{h,j}^k$ given in local coordinates by $a_{h,j}^{k,\ga,\ep}$ (for $\ep \in \{0,1\}^q$)
such that
\[
(a_{h,j} \o \pi_k \o \ps_{k})(x) = x^{j \be} a_{h,j}^{k}(x), \quad (\text{for } 1 \le h \le l \text{ and } 2 \le j \le n_h).
\]

Let us consider the $\cC$-family of polynomials $P^k:=P_1^k \cdots P_l^k$, where
\[
P_h^k(x)(z) := z^{n_h} + \sum_{j=2}^{n_h} (-1)^j a_{h,j}^{k}(x) z^{{n_h}-j}.
\]
Let $y_0:=\ps_k^{-1}(x_0) \in V_k$.
There exist $1 \le h \le l$ and $2 \le j \le n_h$ such that $a_{h,j}^{k}(y_0) \ne 0$, and, thus 
(as $a_{h,1}^{k}=0$), not all roots of $P_h^{k}(y_0)$ coincide. 
Therefore, $|\on{S}( P^k(y_0))| < |\on{S}(P(0))|$, 
and, by the induction hypothesis, 
there exists a finite covering $\{\pi_{kl} : V_{kl} \to V_k\}$ of $V_k$ (possibly shrinking $V_k$) of the type described in $(2)$ such that,
for all $l$, the family of polynomials $P^k \o \pi_{kl}$ admits a $\cC$-parameterization $\la_j^{kl}$ (for $1 \le j \le n$) of its roots on $V_{kl}$.
Since the roots of $P^k$ and $P \o \pi_k \o \ps_k$ differ by the monomial factor $m(x):=x^\be$, 
the $\cC$-functions $x \mapsto m(\pi_{kl}(x)) \cdot \la_j^{kl}(x)$ form a choice of the roots 
of the family $x \mapsto (P \o \pi_k \o \ps_k \o \pi_{kl})(x)$ for $x \in V_{kl}$.
%Then the $\cC$-functions $x \mapsto x^{\be} \la_j^{kl}(x)$ form a choice of the roots 
%of the family $x \mapsto (P \o \pi_k \o \ps_k \o \pi_{kl})(x)$ for $x \in V_{kl}$.

Since $k$ and $x_0$ were arbitrary, the assertion of the theorem follows (by \ref{not}).
\qed\end{demo}

\subsection{Hyperbolic version} \label{hypcCperturb}
If $P$ is hyperbolic, no local power substitutions are needed, see theorem \ref{hypcCperturbth}.

\begin{lemma} \label{hypcCperturblem}
Let $U \subseteq \R^q$ be an open neighborhood of $0$.
Consider a family of hyperbolic polynomials 
\[
P(x)(z) = z^n + \sum_{j=2}^n (-1)^j a_j(x) z^{n-j}, 
\]
with coefficients $a_j$ (for $2 \le j \le n$) in $\cC(U,\R)$.
Assume that $a_2 \ne 0$ and that, for all $j$, $a_j\ne 0$ implies $a_j(x)=x^{\al_j} b_j(x)$, 
where $b_j \in \cC(U,\R)$ is non-vanishing, and $\al_j \in \N^q$.
Then there exists a $\de \in \N^q$ such that $\al_2 = 2 \de$ and $\al_j \ge j \de$, for those $j$ with $a_j \ne 0$.
\end{lemma}

\begin{demo}{Proof}
Since $0 \le \tilde \De_2 (P) = -2 n a_2$ (by theorem \ref{hypth}), we have $\al_2 = 2 \de$ for some $\de \in \N^q$.
If $\de =0$, the assertion is trivial. Let as assume that $\de \ne 0$.

Set $\mu = (\mu_1,\ldots,\mu_q)$, where
\begin{equation} \label{hypmu}
\mu_i := \min\Big\{\frac{(\al_j)_i}{j} : a_j \ne 0\Big\}.
\end{equation}
For contradiction, assume that there is an $i_0$ such that $\mu_{i_0} < \de_{i_0}$.
Consider 
\[
\tilde P(x)(z) := z^n + \sum_{j=2}^n (-1)^j  x^{-j\mu} a_j(x) z^{n-j}.
\]
If all $x_i \ge 0$, then $\tilde P$ is continuous (by \eqref{hypmu}), and if all $x_i > 0$, 
then $\tilde P$ is hyperbolic (its roots differ from those of $P$ by the factor $x^{-\mu}$).
Since the space of hyperbolic polynomials of fixed degree  is closed (by
theorem \ref{hypth}), $\tilde P$ is hyperbolic, if all $x_i \ge 0$.
Since $(\al_2)_{i_0}-2\mu_{i_0}=2\de_{i_0}-2\mu_{i_0}>0$,
all roots (and thus all coefficients) of $\tilde P(x)$ vanish on $\{x_{i_0}=0\}$ 
(as the first and second coefficient vanish, see \ref{hyp}). 
This is a contradiction for those $j$ with $(\al_j)_{i_0}=j \mu_{i_0}$. 
\qed\end{demo}

\begin{theorem}[$\cC$-perturbation of hyperbolic polynomials] \label{hypcCperturbth}
Let $M$ be a $\cC$-mani\-fold. 
Consider a family of hyperbolic polynomials 
\[
P(x)(z) = z^n + \sum_{j=1}^n (-1)^j a_j(x) z^{n-j}, 
\]
with coefficients $a_j$ (for $1 \le j \le n$) in $\cC(M,\R)$.
Let $K \subseteq M$ be compact.
Then there exist:
\begin{enumerate}
\item a neighborhood $W$ of $K$, and
\item a finite covering $\{\pi_k : U_k \to W\}$ of $W$, where each $\pi_k$ 
is a composite of finitely many local blow-ups with smooth center,
\end{enumerate}
such that, for all $k$, the family of polynomials 
$P \o \pi_k$ allows a $\cC$-parameterization of its roots on $U_k$. 
\end{theorem}

\begin{demo}{Proof}
It suffices to modify the proof in \ref{cCperturb} such that no local power substitution is needed.
Suppose we have reduced the problem in virtue of \ref{reduce}.

So $a_{h,j}(0)=0$ for all $1 \le h \le l$ and $1 \le j \le n_h$.
Since $a_{h,1}=0$, we can assume that $a_{h,2} \ne 0$ for all $h$ (otherwise all roots of $P_h$ are identically $0$, see \ref{hyp}).
By theorem \ref{resth}, we find a finite covering $\{\pi_k : U_k \to U\}$ of a neighborhood $U$ of $0$ by $\cC$-mappings $\pi_k$, 
each of which is a composite of finitely many local blow-ups with smooth center, such that, for each $k$, 
the non-zero $a_{h,j} \o \pi_k$ (for $1 \le h \le l$ and $2 \le j \le n_h$)  
simultaneously have only normal crossings.

Let $k$ be fixed and let $x_0 \in U_k$. 
Then $x_0$ admits a neighborhood $W_k$ with suitable coordinates in which $x_0=0$ and such that 
(for $1 \le h \le l$ and $2 \le j \le n_h$) 
either $a_{h,j} \o \pi_k=0$ or 
\begin{equation} \label{hypa}
(a_{h,j} \o \pi_k)(x)=x^{\al_{h,j}} a_{h,j}^{k}(x), 
\end{equation}
where $a_{h,j}^{k}$ is a non-vanishing $\cC$-function on $W_k$, and $\al_{h,j} \in \N^q$.
By lemma \ref{hypcCperturblem}, for each $h$, there exists a $\de_h \in \N^q$
such that $\al_{h,2} = 2 \de_h$.

If some $\de_h =0$, then $(a_{h,2} \o \pi_k)(x_0) = a_{h,2}^{k}(x_0) \ne 0$ and so not all roots of $(P_h \o \pi_k)(x_0)$ coincide.
Thus, $|\on{S}((P\o\pi_k)(x_0))| < |\on{S}(P(0))|$, and, by the induction hypothesis, 
there exists a finite covering $\{\pi_{kl} : W_{kl} \to W_k\}$ of $W_k$ (possibly shrinking $W_k$) of the type described in $(2)$ such that, 
for all $l$, the family of polynomials $P \o \pi_k \o \pi_{kl}$ allows a $\cC$-parameterization of its roots on $W_{kl}$. 

Let us assume that $\de_h \ne 0$ for all $1 \le h \le l$.
By lemma \ref{hypcCperturblem}, we have $\al_{h,j} \ge j \de_h$, for all $1 \le h \le l$ and those $2 \le j \le n_h$ with $a_{h,j} \o \pi_k \ne 0$.
Then 
\[
P_h^{k}(x)(z) := z^{n_h} + \sum_{j=2}^{n_h} (-1)^j x^{-j \de_h} a_{h,j}(\pi_k(x)) z^{{n_h}-j}
\]
is a $\cC$-family of hyperbolic polynomials.
Since $\al_{h,2}=2\de_h$ and $a_{h,2}^k(x_0)\ne 0$, not all roots of $P_h^k(x_0)$ coincide.
Put $P^k := P_1^k \cdots P_l^k$.
Then, $|\on{S}( P^k(x_0))| < |\on{S}(P(0))|$, 
and, by the induction hypothesis, 
there exists a finite covering $\{\pi_{kl} : W_{kl} \to W_k\}$ of $W_k$ (possibly shrinking $W_k$) of the type described in $(2)$ such that,
for all $l$, the family of polynomials $P^k \o \pi_{kl}$ admits a $\cC$-parameterization 
$\la_{h,j}^{kl}$ (for $1 \le h \le l$ and $1 \le j \le n_h$) of its roots on $W_{kl}$.
Since the roots of $P_h^k$ and $P_h \o \pi_k$ differ by the monomial factor $m_h(x):=x^{\de_h}$, 
the $\cC$-functions $x \mapsto m_h(\pi_{kl}(x)) \cdot \la_{h,j}^{kl}(x)$ form a choice of the roots 
of the family $x \mapsto (P \o \pi_k \o \pi_{kl})(x)$ for $x \in W_{kl}$.
%Then the $\cC$-functions $x \mapsto x^{\de_h} \la_{h,j}^{kl}(x)$ form a choice of the roots 
%of the family $x \mapsto (P \o \pi_k \o \pi_{kl})(x)$ for $x \in W_{kl}$.

Since $k$ and $x_0$ were arbitrary, the assertion of the theorem follows (by \ref{not}).
\qed\end{demo}

If the parameter space is one dimensional, we obtain a $\cC$-version of Rellich's classical theorem \cite[Hilfssatz 2]{Rellich37I} 
(see also \cite[5.1]{AKLM98}):

\begin{corollary} \label{hypcCperturbcor}
Let $I \subseteq \R$ be an open interval. Consider a curve of hyperbolic polynomials 
\[
P(x)(z) = z^n + \sum_{j=1}^n (-1)^j a_j(x) z^{n-j}, 
\]
with coefficients $a_j$ (for $1 \le j \le n$) in $\cC(I,\R)$. 
Then there exists a global parameterization $\la_j \in \cC(I,\R)$ (for $1 \le j \le n$) of the roots of $P$.
\end{corollary}

\begin{demo}{Proof}
The local statement follows immediately from theorem \ref{hypcCperturbth}.
(Each local blow-up is the identity map, and, in fact, each non-zero $a_j$ automatically has only normal crossings.)
We claim that a local choice of $\cC$-roots is unique up to permutations. 
In view of this uniqueness property we may glue the local parameterizations of the roots of $P$ to a global one.

For the proof of the claim let $\la^i = (\la_1^i,\ldots,\la_n^i) : J \to \R^n$ (for $i=1,2$) be two local $\cC$-parameterizations 
of the roots of $P$. Let $x_k \to x_\infty \in J$ be a sequence converging in $J$. 
For each $k$ there exists a permutation $\ta_k \in \on{S}_n$ such that $\la^1(x_k) = \ta_k (\la^2(x_k))$. 
Passing to a subsequence, we may assume  that $\la^1(x_k) = \ta (\la^2(x_k))$ for all $k$ and a fixed $\ta \in \on{S}_n$.
By Rolle's theorem (applied repeatedly), the Taylor series at $x_\infty$ of $\la^1$ and $\ta \o \la^2$ coincide.
Quasianalyticity (\ref{cC}.4) implies that $\la^1 = \ta \o \la^2$. 
\qed\end{demo}

\subsection{Real analytic perturbation of polynomials} \label{Comega}
If $\cC=C^\om$, theorem \ref{cCperturb} can be strengthened. 

\begin{theorem}[$C^\om$-perturbation of polynomials] \label{Comegath} 
Let $M$ be a real analytic manifold. 
Consider a family of polynomials 
\[
P(x)(z) = z^n + \sum_{j=1}^n (-1)^j a_j(x) z^{n-j}, 
\]
with coefficients $a_j$ (for $1 \le j \le n$) in $C^\om(M,\C)$.
Let $K \subseteq M$ be compact.
Then there exist:
\begin{enumerate}
\item a neighborhood $W$ of $K$,
\item a finite covering $\{\pi_k : U_k \to W\}$ of $W$, where each $\pi_k$ is a composite of finitely many local blow-ups with smooth center,
\item a finite covering $\{\pi_{kl} : U_{kl} \to U_k\}$ of each $U_k$, where each $\pi_{kl}$ is a single local power substitution.
\end{enumerate}
such that, for all $k,l$, the family of polynomials 
$P \o \pi_k \o \pi_{kl}$ allows a real analytic parameterization of its roots on $U_{kl}$.
\end{theorem}

\begin{demo}{Proof}
Applying resolution of singularities (e.g.\ Hironaka's classical theorem \cite{Hironaka64}, or theorem \ref{resth} for $\cC=C^\om$), 
we obtain that $\tilde \De_s(P \o \pi_k)$ has only normal crossings, 
where $s$ is maximal with the property that $\tilde \De_s(P) \ne 0$ (locally). 
Note that $\tilde \De_s (P)$ is up to a constant factor the discriminant of the square-free reduction of $P$.
Then the assertion follows from the Abhyankar--Jung theorem \cite{Jung08}, \cite{Abhyankar55} 
(see also \cite{KiyekVicente04}, \cite[Section 5]{Sussmann90}, and \cite[Lemma 2.8]{Parusinski94}). 
Here we used that the square-free reduction of a real analytic family of polynomials is real analytic again (see \cite[5.1]{KurdykaPaunescu08}).
\qed\end{demo}

\begin{remarks}
(1) Note that the hyperbolic version of this theorem, where no local power substitutions are needed,
is due to Kurdyka and Paunescu \cite[5.8]{KurdykaPaunescu08}. 

(2) It is unclear to me how to prove this stronger version of theorem \ref{cCperturb} for arbitrary $\cC$ (satisfying (\ref{cC}.1)--(\ref{cC}.6)).
It seems that one can produce a proof of a $\cC$-version of the Abhyankar--Jung theorem along the lines of Luengo's approach \cite{Luengo83}.
Unfortunately, the proof in \cite{Luengo83} contains a gap as pointed out by Kiyek and Vicente \cite{KiyekVicente04}.

(3) Compare this theorem with Parusinski's preparation theorem for subanalytic functions \cite[7.5]{Parusinski94L}.
\end{remarks}

\section{Roots with gradients in \texorpdfstring{$L^1_{\on{loc}}$}{}} \label{secW11}

Let $M$ be a $\cC$-manifold of dimension $q$ equipped with a $C^\infty$ Riemannian metric. 
Consider a family of polynomials
\[
P(x)(z) = z^n + \sum_{j=1}^n (-1)^j a_j(x) z^{n-j},
\]
with coefficients $a_j$ (for $1 \le j \le n$) in $\cC(M,\C)$.
We show in this section that the roots of $P$
admit a parameterization by ``piecewise Sobolev $W^{1,1}_{\on{loc}}$'' functions $\la_i$ (for $1 \le i \le n$). 
That means, there exists a closed nullset $E \subseteq M$ of finite $(q-1)$-dimensional Hausdorff measure such that 
each $\la_i$ belongs to $W^{1,1}(K \setminus E)$ for all compact subsets $K \subseteq M$.
In particular, the classical derivative $\nabla \la_i$ exists almost everywhere and belongs to $L^1_{\on{loc}}$.
The distributional derivatives of the $\la_i$ may not be locally integrable. In fact,
$P$ does in general not allow roots in $W^{1,1}_{\on{loc}}$ (by example \ref{sqroot}). 

\subsection{} \!\!\label{Hmeasure}
We denote by $\cH^k$ the $k$-dimensional Hausdorff measure. It depends on the metric but not on the ambient space.
Recall that for a Lipschitz mapping $f : U \to \R^p$, $U \subseteq \R^q$, we have 
\begin{equation} \label{HmeasureEq}
\cH^k(f(E)) \le \big(\on{Lip}(f)\big)^k \cH^k(E), \quad \text{ for all } E \subseteq U,
\end{equation}
where $\on{Lip}(f)$ denotes the Lipschitz constant of $f$.
The $q$-dimensional Hausdorff measure $\cH^q$ and the $q$-dimensional Lebesgue measure $\cL^q$ coincide in $\R^q$.
If $B$ is a subset of a $k$-plane in $\R^q$ then $\cH^k(B)=\cL^k(B)$.

\subsection{The class \texorpdfstring{$\cW^{\cC}$}{}} \label{classW}
Let $M$ be a $\cC$-manifold of dimension $q$
equipped with a $C^\infty$ Riemannian metric $g$.
We denote by $\cW^{\cC}(M)$ the class of all real or complex valued
functions $f$ with the following properties:
\begin{enumerate}
\item[($\cW_1$)] $f$ is defined and of class $\cC$ on the complement $M \setminus E_{M,f}$ of a closed set $E_{M,f}$ 
with $\cH^q(E_{M,f})=0$ and $\cH^{q-1}(E_{M,f})<\infty$.
\item[($\cW_2$)] $f$ is bounded on $M \setminus E_{M,f}$.
\item[($\cW_3$)] $\nabla f$ belongs to $L^1(M \setminus E_{M,f})=L^1(M)$.
\end{enumerate}

For example, the Heaviside function belongs to $\cW^\cC((-1,1))$, but the function $f(x):= \sin 1/|x|$ does not. 
A $\cW^\cC$-function $f$ may or may not be defined on $E_{M,f}$.
Note that, if the volume of $M$ is finite, then
\begin{equation} \label{Wincl}
f \in \cW^{\cC}(M) \Longrightarrow f \in L^\infty(M \setminus E_{M,f}) \cap W^{1,1}(M \setminus E_{M,f}).
\end{equation}
We shall also use the notations $\cW^{\cC}_{\on{loc}}(M)$ and $\cW^{\cC}(M,\C^n)=(\cW^{\cC}(M,\C))^n$ with the obvious meanings.

In general $\cW^\cC(M)$ depends on the Riemannian metric $g$. 
It is easy to see that $\cW^\cC(U)$ is independent of $g$ for any relatively compact open subset $U \subseteq M$. 
Thus also $\cW^{\cC}_{\on{loc}}(M)$ is independent of $g$.
If $(U,u)$ is a relatively compact coordinate chart and $g_{ij}^u$ is the coordinate expression of $g$, then there exists a constant $C$ 
such that $(1/C) \de_{ij} \le g_{ij}^u \le C \de_{ij}$ as bilinear forms. 

\textbf{From now on, given a $\cC$-manifold $M$, we tacitly choose a $C^\infty$ Riemannian metric $g$ on $M$ and consider $\cW^\cC(M)$ 
with respect to $g$.}

\subsection{} \label{invps}
Let $\rh = (\rh_1,\ldots,\rh_q) \in (\R_{>0})^q$, $\ga=(\ga_1,\ldots,\ga_q) \in (\N_{>0})^q$, 
and $\ep = (\ep_1,\ldots,\ep_q) \in \{0,1\}^q$.
Set
\begin{align*}
\Om(\rh) &:= \{x \in \R^q : |x_j| < \rh_j ~\text{for all}~j\},\\
\Om_\ep(\rh) &:= \{x \in \R^q : 0 < (-1)^{\ep_j} x_j < \rh_j ~\text{for all}~j\}.
\end{align*}
Then $\Om(\rh) \setminus \{\prod_j x_j = 0\} = \bigsqcup \{\Om_\ep(\rh) : \ep \in \{0,1\}^q\}$.
The power transformation
\[
\ps_{\ga,\ep} : \R^q \to \R^q : (x_1,\ldots,x_q) \mapsto ((-1)^{\ep_1} x_1^{\ga_1},\ldots,(-1)^{\ep_q} x_q^{\ga_q})
\] 
maps $\Om_\mu(\rh)$ onto $\Om_\nu(\rh^\ga)$, 
where $\nu=(\nu_1,\ldots,\nu_q)$ such that $\nu_j \equiv \ep_j + \ga_j \mu_j \mod 2$ for all $j$. 
The range of the $j$-th coordinate behaves differently depending on whether $\ga_j$ is even or odd.
So let us consider
\begin{align*}
\bar \ps_{\ga,\ep} : \Om_\ep(\rh) \to \Om_\ep(\rh^\ga) :
(x_1,\ldots,x_q) \mapsto ((-1)^{\ep_1} |x_1|^{\ga_1},\ldots,(-1)^{\ep_q} |x_q|^{\ga_q})
\end{align*}
and its inverse mapping
\begin{align*}
\bar \ps_{\ga,\ep}^{-1} : \Om_\ep(\rh^\ga) \to \Om_\ep(\rh) :
(x_1,\ldots,x_q) \mapsto 
((-1)^{\ep_1} |x_1|^{\frac{1}{\ga_1}},\ldots,(-1)^{\ep_q} |x_q|^{\frac{1}{\ga_q}}).
\end{align*}
Then we have $\bar \ps_{\ga,\ep} \o \bar \ps_{\ga,\ep}^{-1} = \on{id}_{\Om_\ep(\rh^\ga)}$ 
and $\bar \ps_{\ga,\ep}^{-1} \o \bar \ps_{\ga,\ep} = \on{id}_{\Om_\ep(\rh)}$
for all $\ga \in (\R_{>0})^q$ and $\ep \in \{0,1\}^q$.
Note that 
\begin{equation} \label{inclps}
\{\bar \ps_{\ga,\ep} : \ep \in \{0,1\}^q\} \subseteq \{\ps_{\ga,\mu}|_{\Om_\ep(\rh)} : \ep,\mu \in \{0,1\}^q\}.
\end{equation}

\begin{lemma} \label{invpslem}
If $f \in \cW^{\cC}(\Om_\ep(\rh))$ then $f \o \bar \ps_{\ga,\ep}^{-1} \in \cW^{\cC}(\Om_\ep(\rh^\ga))$.
\end{lemma}

\begin{demo}{Proof} 
The mapping $\bar \ps_{\ga,\ep} : \Om_\ep(\rh) \to \Om_\ep(\rh^\ga)$ is a $\cC$-diffeomorphism 
(by (\ref{cC}.1) and (\ref{cC}.6)), and it is Lipschitz. 
Hence, $E_{\Om_\ep(\rh^\ga),f \o \bar \ps_{\ga,\ep}^{-1}}=\bar \ps_{\ga,\ep}(E_{\Om_\ep(\rh),f})$ is closed, 
and we have $\cH^q(E_{\Om_\ep(\rh^\ga),f \o \bar \ps_{\ga,\ep}^{-1}})=0$ and 
$\cH^{q-1}(E_{\Om_\ep(\rh^\ga),f \o \bar \ps_{\ga,\ep}^{-1}}) < \infty$, 
by \eqref{HmeasureEq}. This implies ($\cW_1$) and ($\cW_2$).
Since $f \in \cW^{\cC}(\Om_\ep(\rh))$, we have $\p_i f \in
L^1(\Om_\ep(\rh))$. Thus 
\begin{align*}
\infty &> \int_{\Om_\ep(\rh)} |\p_i f(x)| dx
= \int_{\Om_\ep(\rh^\ga)} |\p_i f(\bar \ps_{\ga,\ep}^{-1}(x))||\det d \bar \ps_{\ga,\ep}^{-1}(x)| dx \\
&=  \Big(\prod_{j \ne i} \tfrac{1}{\ga_j}\Big) 
\int_{\Om_\ep(\rh^\ga)} |\p_i (f \o \bar \ps_{\ga,\ep}^{-1})(x)| \prod_{j \ne i} |x_j|^{\frac{1-\ga_j}{\ga_j}} d x \\
&\ge  \Big(\prod_{j \ne i} \tfrac{\rh_j^{1-\ga_j}}{\ga_j}\Big) \int_{\Om_\ep(\rh^\ga)} |\p_i (f \o \bar \ps_{\ga,\ep}^{-1})(x)| dx.
\end{align*}
That shows ($\cW_3$).
\qed\end{demo}

\subsection{}\label{invP}
Let us define $\bar \ps_\ga^{-1} : \Om(\rh^\ga) \to \Om(\rh)$ by setting $\bar \ps_\ga^{-1}|_{\Om_\ep(\rh^\ga)} := \bar \ps_{\ga,\ep}^{-1}$, 
for $\ep \in \{0,1\}^q$, and by extending it continuously to $\Om(\rh^\ga)$.
Analogously, define $\bar \ps_\ga : \Om(\rh) \to \Om(\rh^\ga)$ such that 
$\bar \ps_{\ga} \o \bar \ps_{\ga}^{-1} = \on{id}_{\Om(\rh^\ga)}$ 
and $\bar \ps_{\ga}^{-1} \o \bar \ps_{\ga} = \on{id}_{\Om(\rh)}$.

Lemma \ref{invpslem} implies:

\begin{lemma} \label{invPlem}
If $f \in \cW^{\cC}(\Om(\rh))$ then $f \o \bar \ps_{\ga}^{-1} \in \cW^{\cC}(\Om(\rh^\ga))$. 
\end{lemma}

\begin{demo}{Proof}
The set 
\[
E_{\Om(\rh^\ga),f \o \bar \ps_{\ga}^{-1}} = \bigcup_{\ep \in \{0,1\}^q} E_{\Om_\ep(\rh^\ga),f \o \bar \ps_{\ga,\ep}^{-1}}  \cup \{x \in \Om(\rh^\ga) : \prod_j x_j=0\}
\] 
obviously has the required properties.
\qed\end{demo}

\subsection{} \label{invvh}
Let $I$ be a subset of $\{1,\ldots,q\}$ with $|I| \ge 2$.
For $i \in I$ consider the mapping $\vh_i : \R^q \to \R^q : x \mapsto y$ given by  
\begin{align} \label{bl}
y_j = \left\{
\begin{array}{l@{\quad \text{for} ~}l}
x_i, & j=i \\
x_i x_j, & j \in I \setminus\{i\} \\
x_j, & j \not\in I
\end{array}
\right..
\end{align}
The image $\vh_i(\Om(\rh) \setminus \{x_i=0\}) =: \tilde \Om_i(\rh)$ has the form
\[
\tilde \Om_i(\rh) = \{x \in \R^q : 0 < |x_i| < \rh_i, |x_j| < \rh_j |x_i| ~\text{for}~j \in I \setminus \{i\}, |x_j| < \rh_j ~\text{for}~j \not\in I\}.
\]
If $\rh_i > 1$ for all $i \in I$, 
then $\Om(\rh) \setminus \{x_i=0 ~\text{for all}~ i \in I\} \subseteq \bigcup_{i \in I} \tilde \Om_i(\rh)$.
Let us consider $\tilde \vh_i := \vh_i|_{\Om(\rh) \setminus \{x_i=0\}}$ and
its inverse mapping
$\tilde \vh_i^{-1} : \tilde \Om_i(\rh) \to \Om(\rh) \setminus \{x_i=0\} : x \mapsto y$ given by
\begin{align*} 
y_j = \left\{
\begin{array}{l@{\quad \text{for} ~}l}
x_i, & j=i \\
\frac{x_j}{x_i}, & j \in I \setminus\{i\} \\
x_j, & j \not\in I
\end{array}
\right..
\end{align*}

\begin{lemma} \label{invvhlem} 
If $f \in \cW^{\cC}(\Om(\rh))$ then $f \o \tilde \vh_i^{-1} \in \cW^{\cC}(\tilde \Om_i(\rh))$.
\end{lemma}

\begin{demo}{Proof}
We may view $f$ as a function in $\cW^\cC(\Om(\rh) \setminus \{x_i=0\})$, 
where $E_{\Om(\rh) \setminus \{x_i=0\},f}= E_{\Om(\rh),f} \setminus \{x_i =0\}$.
The mapping $\tilde \vh_i : \Om(\rh) \setminus \{x_i=0\} \to \tilde \Om_i(\rh)$ is a $\cC$-diffeomorphism 
(by (\ref{cC}.1) and (\ref{cC}.6)), and it is Lipschitz. 
Hence, $E_{\tilde \Om_i(\rh),f \o \tilde \vh_i^{-1}}=\tilde \vh_i(E_{\Om(\rh)\setminus \{x_i=0\},f})$ is closed, 
and we have $\cH^q(E_{\tilde \Om_i(\rh),f \o \tilde\vh_i^{-1}})=0$ and $\cH^{q-1}(E_{\tilde \Om_i(\rh),f \o \tilde \vh_i^{-1}}) < \infty$, 
by \eqref{HmeasureEq}. This implies ($\cW_1$) and ($\cW_2$).

The following identities are consequences of the substitution formula (applied from right to left).
The right-hand sides are finite, since $\p_j f \in L^1(\Om(\rh))$ for all $j$ and since $|I| \ge 2$.
\begin{align*}
&\int_{\tilde \Om_i(\rh)} \left|\p_i f(\tilde \vh_i^{-1}(x))\right| dx 
= \int_{\Om(\rh)} |\p_i f(x)| |x_i|^{|I|-1} dx < \infty, \\
&\int_{\tilde \Om_i(\rh)} \left|\p_j f(\tilde \vh_i^{-1}(x)) \frac{x_j}{x_i^2}\right| dx 
= \int_{\Om(\rh)} |\p_j f(x)| |x_i|^{|I|-2} |x_j| dx < \infty, \quad j \in I \setminus \{i\}, \\
&\int_{\tilde \Om_i(\rh)} \left|\p_j f(\tilde \vh_i^{-1}(x)) \frac{1}{x_i}\right| dx
= \int_{\Om(\rh)} |\p_j f(x)| |x_i|^{|I|-2} dx < \infty, \quad j \in I \setminus \{i\}, \\
&\int_{\tilde \Om_i(\rh)} \left|\p_j f(\tilde \vh_i^{-1}(x))\right| dx
= \int_{\Om(\rh)} |\p_j f(x)| |x_i|^{|I|-1} dx < \infty, \quad j \not\in I. 
\end{align*}
It follows that the partial derivatives
\begin{align*}
\p_j (f \o \tilde \vh_i^{-1})(x) = \left\{ 
\begin{array}{l@{\quad \text{for} ~}l}
\p_i f(\tilde \vh_i^{-1}(x)) -\sum_{k\in I \setminus\{i\}} \p_k f(\tilde \vh_i^{-1}(x)) \frac{x_k}{x_i^2}, & j=i\\
\p_j f(\tilde \vh_i^{-1}(x)) \frac{1}{x_i}, & j \in I \setminus \{i\} \\
\p_j f(\tilde \vh_i^{-1}(x)), & j \not\in I
\end{array}
\right.
\end{align*}
belong to $L^1(\tilde \Om_i(r))$. Thus $(\cW_3)$ is shown.
\qed\end{demo}

\begin{lemma} \label{invbl}
Let $\vh : M' \to M$ be a blow-up of a $\cC$-manifold $M$ with center a closed $\cC$-submanifold $C$ of $M$.
If $f \in \cW^{\cC}_{\on{loc}}(M')$ then $f \o (\vh|_{M' \setminus \vh^{-1}(C)})^{-1} \in 
\cW^{\cC}_{\on{loc}}(M)$. 
\end{lemma}

\begin{demo}{Proof}
Let $K \subseteq M$ be compact. Hence $K$ can be covered by finitely many relatively compact coordinate neighborhoods $(U,u)$ 
such that over $U$ the mapping 
$\vh$ identifies with the mapping $U' \to U$ described in \ref{blowup}.
Each $U'$ is covered by charts $(U_i',u_i')$ such that $u \o \vh|_{U_i'} \o (u_i')^{-1} = \vh_i$ (where $\vh_i$ is defined in \eqref{bl}).

\[
\xymatrix{
M' \ar@{->>}[d]_{\vh} & ~U' \ar@{_(->}[l] \ar@{->>}[d]_{\vh|_{U'}} & ~U_i' \ar@{_(->}[l] \ar[r]^{u_i'} \ar@{->>}[d]_{\vh|_{U_i'}} 
& \Om(\rh) \ar@{->>}[d]_{\vh_i} & ~\Om(\rh)\setminus \{x_i=0\} \ar@{_(->}[l] \ar@{->>}[d]_{\tilde \vh_i}\\
M & ~U \ar@{_(->}[l] & ~\vh(U_i') \ar@{_(->}[l] \ar[r]_{u|_{\vh(U_i')}} & \vh_i(\Om(\rh)) & ~\tilde \Om_i(\rh) \ar@{_(->}[l]
}
\]

Since $\vh$ is proper and $U$ is relatively compact, $U'$ is relatively compact as well.
Thus $f|_{U'} \in \cW^\cC(U')$, and $\cW^\cC(U')$ is independent
of the Riemannian metric. 
We may assume that there is a $\rh \in (\R_{>1})^q$ such that $u_i'(U_i') = \Om(\rh)$.
By lemma \ref{invvhlem},  $f|_{U'_i} \o (u_i')^{-1} \o \tilde \vh_i^{-1} \in \cW^{\cC}(\tilde \Om_i(\rh))$.
Since $u_i'(U_i' \setminus \vh^{-1}(C)) = \Om(\rh) \setminus \{x_i=0\}$ and 
$\tilde \vh_i=\vh_i|_{\Om(\rh) \setminus \{x_i=0\}}$, we have
\begin{equation} \label{f_i}
f|_{U'_i} \o (u_i')^{-1} \o \tilde \vh_i^{-1} = f|_{U'_i} \o (\vh|_{U_i' \setminus \vh^{-1}(C)})^{-1} \o u^{-1}|_{\tilde \Om_i(\rh)}
\in \cW^{\cC}(\tilde \Om_i(\rh)).
\end{equation}
Let $\Upsilon(\rh) := \bigcup_{i \in I} \tilde \Om_i(\rh)$. Note that $\Om(\rh) \setminus \{x_i=0 ~\text{for all}~ i \in I\} \subseteq \Upsilon(\rh)$.
Then
\begin{equation} \label{f}
f|_{U'} \o (\vh|_{U' \setminus \vh^{-1}(C)})^{-1} \o u^{-1}|_{\Upsilon(\rh)}
\in \cW^{\cC}(\Upsilon(\rh)),
\end{equation}
where $E_{\Upsilon(\rh),\star} := \bigcup_{i \in I} \big(E_{\tilde \Om_i(\rh),\star\star} \cup \p \big(\tilde \Om_i(\rh)\big)\big)$
and $\star$ and $\star\star$ represent the functions in \eqref{f} and \eqref{f_i}, respectively.
So we find (possibly shrinking $U$)
\[
f \o (\vh|_{M' \setminus \vh^{-1}(C)})^{-1}|_U=
f|_{U'} \o (\vh|_{U' \setminus \vh^{-1}(C)})^{-1} 
\in \cW^{\cC}(U),
\]
where $\cW^{\cC}(U)$ is independent of the Riemannian metric.
It follows immediately that 
\[
f \o (\vh|_{M' \setminus \vh^{-1}(C)})^{-1}|_{\bigcup U}  \in \cW^{\cC}(\bigcup U),  
\]
where the union in finite.
This completes the proof.
\qed\end{demo}

\begin{lemma} \label{patch}
Let $K \subseteq M$ be compact, let $\{(U_j,u_j) : 1 \le j \le N\}$ be a finite collection of connected relatively compact coordinate charts covering $K$, and
let $f_j \in \cW^{\cC}(U_j)$.
Then, after shrinking the $U_j$ slightly such that they still cover $K$, there exists a function $f \in \cW^\cC(\bigcup_j U_j)$ satisfying the following
condition:
\begin{enumerate}
\item[(1)] If $x \in \bigcup_j U_j$ then either $x \in E_{\bigcup_j U_j}$ or $f(x)=f_j(x)$ for a $j \in \{i : x \in U_i\}$.
\end{enumerate}
\end{lemma}

\begin{demo}{Proof}
We construct $f$ step-by-step.
Suppose that a function $f' \in \cW^\cC(\bigcup_{j=1}^{n-1} U_j)$ satisfying (1) has been found.
If $(\bigcup_{j=1}^{n-1} U_j) \cap U_n = \emptyset$ then the function 
\[
f := f' \mathbf{1}_{\bigcup_{j=1}^{n-1} U_j} + f_n \mathbf1_{U_n} \in \cW^\cC(\bigcup\nolimits_{j=1}^{n} U_j)
\] 
has property (1). 
Otherwise, consider the chart $(U_n,u_n)$. We may assume that $u_n(U_n) = B_1(0)$, the open unit ball in $\R^q$. Choose $\ep>0$ small, 
such that the collection $\{U_j : 1 \le j \le N, j \ne n\} \cup U_n'$, where $U_n' := u_n^{-1}(B_{1-\ep}(0))$, still covers $K$.
The set $S := \p B_{1-\ep}(0) \cap u_n((\bigcup_{j=1}^{n-1} U_j) \cap U_n)$ is closed in $u_n((\bigcup_{j=1}^{n-1} U_j) \cap U_n)$,
$\cH^q(S)=0$, and $\cH^{q-1}(S)< \infty$. 
So $u_n^{-1}(S)$ is closed in $\bigcup_{j=1}^{n-1} U_j \cup U_n'$, and, by \eqref{HmeasureEq}, 
$\cH^q(u_n^{-1}(S))=0$, and $\cH^{q-1}(u_n^{-1}(S))< \infty$.
Thus 
\[
f := f' \mathbf{1}_{(\bigcup_{j=1}^{n-1} U_j) \setminus U_n'} + f_n \mathbf1_{U_n'} \in \cW^\cC(\bigcup\nolimits_{j=1}^{n-1} U_j \cup U_n')
\]
and satisfies (1). Repeating this procedure finitely many times, produces the required function.
\qed\end{demo}

\begin{theorem}[$\cW^{\cC}$-roots] \label{Wloc} 
Let $M$ be a $\cC$-manifold.
Consider a family of polynomials
\[
P(x)(z) = z^n + \sum_{j=1}^n (-1)^j a_j(x) z^{n-j},
\]
with coefficients $a_j$ (for $1 \le j \le n$) in $\cC(M,\C)$.
For any compact subset $K \subseteq M$ there exists a relatively compact neighborhood $W$ of $K$ and a parameterization 
$\la_j$ (for $1 \le j \le n$) of the roots of $P$ on $W$ such that $\la_j \in \cW^{\cC}(W)$ for all $j$.  
In particular, for each $\la_j$ we have $\nabla \la_j \in L^1(W)$.
\end{theorem}

\begin{demo}{Proof}
By theorem \ref{cCperturb}, 
there exists a neighborhood $W$ of $K$ and
a finite covering $\{\pi_k : U_k \to W\}$ of $W$, where each $\pi_k$ is a 
composite of finitely many mappings each of which is either a local blow-up $\Ph$ with smooth 
center or a local power substitution $\Ps$ (cf.\ \ref{not}),
such that, for all $k$, the family of polynomials 
$P \o \pi_k$ allows a $\cC$-parameterization $\la_i^k$ (for $1 \le i \le n$) of its roots on $U_k$.

In view of lemma \ref{patch}, the proof of the theorem will be complete once the following assertions are shown:
\begin{enumerate}
\item Let $\Ps = \iota \o \ps : V \to W \to M$ be a local power substitution. 
If the roots of $P \o \Ps$ allow a parameterization in $\cW^{\cC}_{\on{loc}}$, then so do the roots of $P|_W$.
\item Let $\Ph = \iota \o \vh : U' \to U \to M$ be local blow-up with smooth center.
If the roots of $P \o \Ph$ allow a parameterization in $\cW^{\cC}_{\on{loc}}$, then so do the roots of $P|_U$.
\end{enumerate}

Assertion $(2)$ is an immediate consequence of lemma \ref{invbl}.
To prove $(1)$, let $\la_i^{\Ps}=\la_i^{\ps_{\ga,\ep}}$ (for some $\ga \in (\N_{>0})^q$ and all $\ep \in \{0,1\}^q$, cf.\ \ref{not}) 
be functions in $\cW^{\cC}_{\on{loc}}(V)$ which parameterize the roots of $P \o \Ps$.
We can assume without loss (possibly shrinking $V$) that $V=\Om(\rh)$, $W = \Om(\rh^{\ga})$, 
and that each $\la_i^{\ps_{\ga,\ep}} \in \cW^{\cC}(\Om(\rh))$, for some $\rh \in (\R_{>0})^q$.
Let us define $\la_i^{\bar \ps_{\ga}} \in \cW^{\cC}(\Om(\rh))$ by setting (in view of \eqref{inclps} and \ref{invP})
\[
\la_i^{\bar \ps_\ga}|_{\Om_\ep(\rh)} := \la_i^{\bar \ps_{\ga,\ep}}|_{\Om_\ep(\rh)}, 
\quad \ep \in \{0,1\}^q.
\]
On the set $\{x \in \Om(\rh) : \prod_j x_j = 0\}$ we may define $\la_i^{\bar \ps_{\ga}}$ (for $1 \le i \le n$) arbitrarily 
such that they form a parameterization of the roots of 
$P \o \iota \o \bar \ps_{\ga}$.
By lemma \ref{invPlem}, 
\[
\la_i := \la_i^{\bar \ps_{\ga}} \o \bar \ps_{\ga}^{-1} \in \cW^{\cC}(\Om(\rh^{\ga})) = \cW^{\cC}(W).
\]
Clearly, $\la_i$ (for $1 \le i \le n$) constitutes a parameterization of the roots of $P|_W$. 
Thus the proof of $(1)$ is complete.
\qed\end{demo}

\begin{corollary}[Local $\cW^{\cC}$-sections] \label{corW}
The mapping $\si^n : \C^n\to \C^n$ from roots to coefficients
(cf.\ \eqref{esf}) admits local $\cW^{\cC}$-sections, 
for $\cC$ any class of $C^\infty$-functions satisfying \emph{(\ref{cC}.1)--(\ref{cC}.6)}.
\end{corollary}

\begin{demo}{Proof}
Apply theorem \ref{Wloc} to the family
\[
P(a)(z) = z^n + \sum_{j=1}^n (-1)^j a_j z^{n-j}, \quad a = (a_1,\ldots,a_n) \in \C^n = \R^{2n}. 
\]
It is a $\cC$-family by (\ref{cC}.1).
\qed\end{demo}

In the following we show that the conclusion of theorem \ref{Wloc} is best possible.

\begin{example}[The derivatives of the roots are not in $L^p_{\on{loc}}$ for any $1<p\le \infty$] 
\label{optimal}
\hfill
In general the roots of a $\cC$ (even polynomial)
family of polynomials $P$ do not allow parameterizations $\la_j$ with
$\nabla \la_j \in L^p_{\on{loc}}$
for any $1 < p \le \infty$.
That is shown by the example 
\[ 
P(x)(z) = z^n- x_1 \cdots x_q, \quad  x=(x_1,\ldots,x_q) \in \R^q,
\] 
if $n \ge \frac{p}{p-1}$, for $1 < p <
\infty$,
and if $n \ge 2$, for $p=\infty$.
\end{example}

\subsection{Remark} \label{Lw}
Compare theorem \ref{Wloc} with the results obtained in \cite{CJS83} and
\cite{CL03}:
For a non-negative real valued function $f \in C^k(U)$, where $U \subseteq
\R^q$ is open and $k \ge 2$,
they find in \cite{CJS83} that $\nabla (f^{1/k}) \in L^1_{\on{loc}}(U)$.
Actually, for each compact $K \subseteq U$, one has $\nabla (f^{1/k}) \in
L^{k/(k-2)}_w(K)$, due to \cite{CL03}, 
where $L^p_w$ denotes the \emph{weak} $L^p$ space.
By example \ref{optimal}, we can in general not expect that the 
derivatives of the roots of $P$ belong to any $L^p_w(K)$ with $p>1$, 
since $L^p(K) \subseteq L^p_w(K) \subseteq L^q(K)$ for $1 \le q < p < \infty$.

\subsection{The one dimensional case} \label{1dim}
Let $P$ be a curve of polynomials. 
Then the proof of lemma \ref{invpslem} actually shows that
pullback by $\bar \ps_{\ga,\ep}^{-1}(x) = (-1)^\ep |x|^{1/\ga}$, ($x \in \R$, $\ga \in \N_{>0}$, and $\ep=0,1$), 
preserves absolute continuity.
So theorem \ref{Wloc} reproduces (for $\cC$-coefficients) the following result proved in \cite{RainerAC} 
(see also \cite{Spagnolo99}):

\begin{theorem}
The roots of an everywhere normally nonflat $C^\infty$-curve of polynomials $P$ 
may be parameterized by locally absolutely continuous functions.
\end{theorem}

A curve of polynomials $P$ with $C^\infty$-coefficients $a_j$ is 
\emph{normally nonflat} at $x_0$ if $x \mapsto \tilde \De_s(P(x))$ is not infinitely flat 
at $x_0$, where $s$ is maximal with the property that the germ at $x_0$ of 
$x \mapsto \tilde \De_s(P(x))$ is not $0$. 
Or, equivalently, no two of the continuously chosen roots (which is always possible in the one dimensional case, cf.\ \cite[II 5.2]{Kato76}) 
meet of infinite order of flatness.

On an interval $I \subseteq \R$ the space of locally absolutely continuous functions
coincides with the
Sobolev space $W_{\on{loc}}^{1,1}(I)$. 
However:

\begin{example}[The roots are not in $W^{1,1}_{\on{loc}}$] \label{sqroot}
Multiparameter $\cC$ (even polynomial) families of polynomials do not 
allow roots in $W_{\on{loc}}^{1,1}$, as the following example shows:
\[
P(x)(z) = z^2-x,\quad x\in \C=\R^2.
\]
The roots are
$\la_{12}=\pm\sqrt{x}$ which must have a jump along some ray. The
distributional derivative of $\sqrt{x}$ with respect to angle contains a
delta distribution which is not in $L^1_{\text{loc}}$.
\end{example}

\begin{example}[The roots are not in $VMO$] \label{VMO}
Let $U \subseteq \R^q$ be open. We say that a real or complex valued $f \in L^1_{\on{loc}}(U)$ has \emph{vanishing mean oscillation}, 
or $f \in VMO(U)$, if, for cubes $Q \subseteq \R^q$ with closure $\overline Q \subseteq U$, we have 
\[
\|f\|_{BMO} := \sup\{\on{mo}(f,Q) : Q\} < \infty \quad \text{and} \quad
\lim_{s \to 0} \sup\{\on{mo}(f,Q): |Q| \le s\} = 0,
\]
where 
\[
f_Q := \frac{1}{|Q|} \int_{Q} f(x) dx \quad \text{and} \quad \on{mo}(f,Q) := \frac{1}{|Q|} \int_{Q} |f(x) - f_Q| dx.
\]
Functions $f \in L^1_{\on{loc}}(U)$ with $\|f\|_{BMO} < \infty$ are said to have \emph{bounded mean oscillation} 
(or $f \in BMO(U)$). Cf.\ \cite{Sarason75} and \cite{BN95, BN96}.

By proposition \ref{bdroot}, the roots of a family of polynomials $P$ whose coefficients are bounded functions on $U$
are bounded as well and hence in $BMO(U)$.
Thus it makes sense to ask whether the roots of a $\cC$-family $P$ admit parameterizations in $VMO$. In general the answer is no:
\ref{sqroot} provides a counter example. 

Namely: Let $S=(-\infty,0] \times \{0\} \subseteq \R^2$ be the left $x$-axis and let $f : \R^2 \setminus S \to \C$ be defined, in polar coordinates 
$(r,\ph) \in (0,\infty) \times (-\pi,\pi)$, by 
\[
f(r,\ph) = \sqrt{r}\big(\cos \tfrac{\ph}{2} + i \sin \tfrac{\ph}{2}\big).
\]
For convenience of computation we use
\[
Q(x_0,\ep) := \{(r,\ph) : |r-x_0| < \ep, -\pi< \ph < -\pi+\ep \text{ or } \pi-\ep < \ph < \pi\}, 
\]
where $0 < \ep < x_0 < \pi/2$. 
Since $Q(x_0,\ep)$ is symmetric with respect to the $x$-axis, we find $\Im f_{Q(x_0,\ep)} = (\Im f)_{Q(x_0,\ep)} =0$. 
It is easy to compute
\[
\on{mo}(\Im f,Q(x_0,\ep)) = \tfrac{2}{5} \sin \tfrac{\ep}{2} \cdot \frac{(x_0+\ep)^{\tfrac{5}{2}}-(x_0-\ep)^{\tfrac{5}{2}}}{x_0 \ep^2} 
\East{\ep \to 0}{} \sqrt{x_0}.
\]
Since $\on{mo}(f,Q(x_0,\ep)) \ge \on{mo}(\Im f,Q(x_0,\ep))$, we may conclude that
$f \not\in VMO(U)$, for each open $U \subseteq \R^2$ containing the origin.
\end{example}

\section{Roots with locally bounded variation} \label{secBV}

The roots of a $\cC$-family of polynomials admit a parameterization by functions having locally bounded variation, actually, 
even by $SBV_{\on{loc}}$-functions.

\subsection{Functions of bounded variation} \label{BV} Cf.\ \cite{AFP00}. 
Let $U \subseteq \R^q$ be open. A real valued function $f \in L^1(U)$ is said to have \emph{bounded variation}, or to belong to $BV(U)$, 
if its distributional derivative is representable by a finite Radon measure in $U$, i.e.,
\[
\int_U f \, \p_i \ph ~dx = - \int_U \ph ~d D_i f, \quad \text{ for all } \ph \in C^\infty_c(U) \text{ and } 1 \le i \le q, 
\]
for some $\R^q$-valued measure $D f = (D_1 f,\cdots,D_q f)$ in $U$.
Then $W^{1,1}(U) \subseteq BV(U)$: for any $f \in W^{1,1}(U)$ the distributional derivative is given by $(\nabla f) \cL^q$.
See \cite[Section 3.1]{AFP00} for equivalent definitions and properties of $BV$-functions.

A complex valued function $f : U \to \C$ is said to be of \emph{bounded variation}, or to be in $BV(U,\C)$, if $(\Re f,\Im f) \in (BV(U))^2$.

\subsection{Special functions of bounded variation} \label{SBV} 
This notion is due to \cite{AmbrosioDeGiorgi88}.
For a detailed treatment see \cite{AFP00}.
Let $U \subseteq \R^q$ be open and let $f \in BV(U)$. We may write %by Radon-Nikodym
\[
D f = D^a f + D^s f,
\]
where $D^a f$ is the absolutely continuous part of $D f$ and $D^s f$ is the singular part of $D f$ with respect to $\cL^q$.

We say that $f$ has an \emph{approximate limit at $x \in U$} if there exists $a \in \R$ such that
\[
\lim_{r\searrow 0} \frac{1}{|B_r(x)|} \int_{B_r(x)} |f(y)-a| dy=0.  
\]
The \emph{approximate discontinuity set} $S_f$ is the set of points where this property does not hold. 
A point $x \in U$ is called \emph{approximate jump point of $f$} if there exist $a^\pm \in \R$ and $\nu \in S^{q-1}$ such that $a^+\ne a^-$ and
\[
\lim_{r\searrow 0} \frac{1}{|B_r^\pm(x,\nu)|} \int_{B_r^\pm(x,\nu)} |f(y)-a^\pm| dy=0,
\]
where $B_r^\pm(x,\nu) := \{y \in B_r(x) : \pm \< y-x \mid \nu \> > 0\}$. The set of approximate jump points is denoted by $J_f$.

For any $f \in BV(U)$ the measures
\[
D^j f :=\mathbf{1}_{J_f} D^s f \quad \text{and} \quad D^c f := \mathbf{1}_{U \setminus S_f} D^s f
\]
are called the \emph{jump part} and the \emph{Cantor part} of the derivative. 
Since $D f$ vanishes on the $\cH^{q-1}$-negligible set $S_f \setminus J_f$, we obtain the decomposition 
\[
D f = D^a f + D^j f + D^c f.
\]

We say that $f \in BV(U)$ is a \emph{special function of bounded variation}, and we write $f \in SBV(U)$, if
the Cantor part of its derivative $D^c f$ is zero.

\begin{proposition}[{\cite[4.4]{AFP00}}] \label{SBVprop}
Let $U \subseteq \R^q$ be open and bounded, $E \subseteq \R^q$ closed, and $\cH^{q-1}(E \cap U) < \infty$. 
Then, any function $f : U \to \R$ that belongs to $L^\infty(U \setminus E) \cap W^{1,1}(U \setminus E)$ belongs also to $SBV(U)$ and satisfies 
$\cH^{q-1}(S_f \setminus E)=0$.
\end{proposition}

A complex valued function $f$ belongs to $SBV(U,\C)$ if $(\on{Re} f,\on{Im} f) \in (SBV(U))^2$.

\begin{theorem}[$SBV$-roots] \label{BVth} 
Let $U \subseteq \R^q$ be open.
Consider a family of polynomials
\[
P(x)(z) = z^n + \sum_{j=1}^n (-1)^j a_j(x) z^{n-j},
\]
with coefficients $a_j$ (for $1 \le j \le n$) in $\cC(U,\C)$.
For any compact subset $K \subseteq U$ there exists a relatively compact neighborhood $W$ of $K$ and a parameterization 
$\la_j$ (for $1 \le j \le n$) of the roots of $P$ on $W$ such that $\la_j \in SBV(W,\C)$ for all $j$.
\end{theorem}

\begin{demo}{Proof}
It follows immediately from theorem \ref{Wloc}, proposition \ref{SBVprop}, and \eqref{Wincl}.
\qed\end{demo}

Combining corollary \ref{corW} with proposition \ref{SBVprop} or applying theorem \ref{BVth} to the family $P$ in \ref{corW} gives:

\begin{corollary}[Local $SBV$-sections]
The mapping $\si^n : \C^n\to \C^n$ from roots to coefficients
(see \eqref{esf}) admits local $SBV$-sections, 
for $\cC$ any class of $C^\infty$-functions satisfying \emph{(\ref{cC}.1)--(\ref{cC}.6)}. \qed
\end{corollary}

\section{Perturbation of normal matrices} \label{secmatrix}

We investigate the consequences of our results in the perturbation theory of normal matrices. 
It is evident that the eigenvalues  of a $\cC$-family of normal matrices possess the regularity properties of the roots of a 
$\cC$-family of polynomials. We prove that the same it true for the eigenvectors.

\begin{theorem}[$\cC$-perturbation of normal matrices] \label{qamatrix}
Let $M$ be a $\cC$-manifold. 
Consider a family of normal complex 
matrices 
\[
A(x)=(A_{ij}(x))_{1 \le i,j \le n}
\]
(acting on a complex vector space $V = \C^n$), where the entries $A_{ij}$ (for $1 \le i,j \le n$) 
belong to $\cC(M,\C)$. Let $K \subseteq M$ be compact. 
Then there exist:
\begin{enumerate}
\item a neighborhood $W$ of $K$, and
\item a finite covering $\{\pi_k : U_k \to W\}$ of $W$, where each  
$\pi_k$ is a composite of finitely many mappings each of which is either a 
local blow-up with smooth center or a local power substitution,
\end{enumerate}
such that, for all $k$, the family of normal complex matrices
$A \o \pi_k$ allows a $\cC$-parameterization of its eigenvalues and eigenvectors.

If $A$ is a family of Hermitian matrices, then the above statement holds with each 
$\pi_k$ being a composite of finitely many local blow-ups with smooth center only.
\end{theorem}

\begin{demo}{Proof}
By theorem \ref{cCperturb} applied to the characteristic polynomial
\begin{align} \label{ch}
\ch(A(x))(\la) = \det (A(x) - \la \I) &= \sum_{j=0}^n (-1)^{n-j} \on{Trace}(\La^j
A(x)) \la^{n-j} \\
&=: (-1)^n \Big(\la^n + \sum_{j=1}^n (-1)^j a_j(x) \la^{n-j}\Big), \nonumber
\end{align}
there exist a neighborhood $W$ of $K$ and a finite covering $\{\pi_k : U_k \to W\}$ of $W$ of the type described in (2)
such that, for all $k$, the family of normal matrices $A \o \pi_k$ 
admits a $\cC$-parameterization $\la_i$ (for $1 \le i \le n$) of its eigenvalues.

Let us prove the statement about the eigenvectors. 
We shall show that (for each $k$) there exists a finite covering $\{\pi_{kl} : U_{kl} \to U_k\}$ of $U_k$ 
of the type described in $(2)$ such that 
$A \o \pi_k \o \pi_{kl}$ admits a $\cC$-parameterization of its eigenvectors (for all $l$). 
This assertion follows from the following claim.
Composing the finite coverings in the sense of \ref{not}, will complete the proof.

\begin{claim*}
Let $A=A(x)$ be a family of normal complex $n \times n$ matrices, where the entries $A_{ij}$ are $\cC$-functions 
and the eigenvalues of $A$ admit a $\cC$-parameterization $\la_j$ in a neighborhood of $0 \in \R^q$.
Then there exists a finite covering $\{\pi_k : U_k \to U\}$ of a neighborhood $U$ of $0$ of the type described in $(2)$ such that, for all $k$, 
$A \o \pi_k$ admits a $\cC$-parameterization of its eigenvectors. 
\end{claim*}

\begin{demo}{Proof of the claim}
We use induction on $|\on{S}(\ch(A(0)))|$ (cf.\ \ref{reduce}).

First consider the following reduction:
Let $\nu_1,\ldots,\nu_l$ denote
the pairwise distinct eigenvalues of $A(0)$
with respective multiplicities $m_1,\ldots,m_l$. 
The sets
\[
\La_h:=\{\la_i : \la_i(0)=\nu_h\}, \quad 1 \le h \le l, 
\]
form a partition of the $\la_i$ such that, for $x$ near $0$, $\la_i(x) \ne \la_j(x)$ if $\la_i$ and $\la_j$ belong to different $\La_h$. 
Consider
\begin{align*}
V_x^{(h)} := \bigoplus_{\la \in \La_h}
\on{ker}(A(x)-\la(x)) 
= \on{ker} \big(\circ_{\la \in \La_h}
(A(x)-\la(x))\big), \quad 1 \le h \le l.
\end{align*}
(The order of the compositions is not
relevant.)
So $V_x^{(h)}$ is the kernel of a vector bundle homomorphism $B(x)$ of class $\cC$
with constant rank (even of constant dimension of the kernel),
and thus it is a vector subbundle of class $\cC$ of the trivial bundle 
$U \times V \to U$ (where $U\subseteq \R^q$ is a neighborhood of $0$) which admits a $\cC$-framing.
This can be seen as follows: Choose a basis of $V$ such that $A(0)$
is diagonal. By the elimination procedure
one can construct a basis for the kernel of $B(0)$. For $x$ near $0$, the
elimination procedure (with the same choices)
gives then a basis of the kernel of $B(x)$. 
This clearly involves only operations which preserve the class $\cC$.
The elements of this basis are
then of class $\cC$ in $x$ near $0$.

Therefore, it suffices to find $\cC$-eigenvectors in each subbundle
$V^{(h)}$ separately, expanded in
the constructed frame field of class $\cC$. But in this frame field the vector
subbundle looks again like a constant vector space.
So we may treat each of these parts ($A$ restricted to $V^{(h)}$, as
matrix with respect to the frame field) separately. 
For simplicity of notation we suppress the index $h$.

Let us suppose that all eigenvalues of $A(0)$ coincide and are equal to $a_1(0)/n$, 
according to \eqref{ch}.
Eigenvectors of $A(x)$ are also eigenvectors of $A(x) - (a_1(x)/n) \I$ (and vice versa),
thus we may replace $A(x)$ by $A(x) - (a_1(x)/n) \I$ and assume that the first
coefficient of
the characteristic polynomial \eqref{ch} vanishes identically. Then
$A(0)=0$.

If $A=0$ identically, we choose the eigenvectors constant and we are done.
Note that this proves the claim, if $|\on{S}(\ch(A(0)))|=1$.

Assume that $A \ne 0$.
By theorem \ref{resth} (and \ref{RtoC}), there exists a finite covering 
$\{\pi_k : U_k \to U\}$ 
of a neighborhood $U$ of $0$ by $\cC$-mappings $\pi_k$, each of which is a 
composite of finitely many local blow-ups with smooth center, 
such that, for each $k$,
the non-zero entries $A_{ij} \o \pi_k$ of $A \o \pi_k$ and its pairwise non-zero differences 
$A_{ij} \o \pi_k - A_{lm} \o \pi_k$ simultaneously have only normal crossings. 

Let $k$ be fixed and let $x_0 \in U_k$.
Then $x_0$ admits a neighborhood $W_k$ with suitable coordinates in which $x_0=0$ and such that 
either $A_{ij} \o \pi_k =0$ or 
\[
(A_{ij} \o \pi_k)(x) = x^{\al_{ij}} B_{ij}^k(x), 
\]
where $B_{ij}^k$ is a non-vanishing $\cC$-function on $W_k$, and $\al_{ij} \in \N^q$.
The collection of multi-indices $\{\al_{ij} : A_{ij} \o \pi_k \ne 0\}$ 
is totally ordered, by lemma \ref{order}. Let $\al$ denote its minimum.

If $\al=0$, then $(A_{ij} \o \pi_k)(x_0)=B_{ij}^k(x_0)\ne 0$ for some $1 \le i,j \le n$. 
Since the first coefficient of $\ch(A \o \pi_k)$ vanishes, we may conclude that not all eigenvalues of 
$(A \o \pi_k)(x_0)$ coincide.
Thus, $|\on{S}(\ch(A \o \pi_k)(x_0))| < |\on{S}(\ch(A(0)))|$, and, by the induction hypothesis, 
there exists a finite covering $\{\pi_{kl} : W_{kl} \to W_k\}$ of $W_k$ (possibly shrinking $W_k$) of the type described in $(2)$ such that, 
for all $l$, the family of normal matrices $A \o \pi_k \o \pi_{kl}$ allows a $\cC$-parameterization of its eigenvectors on $W_{kl}$.

Assume that $\al \ne 0$.
Then there exist $\cC$-functions $A_{ij}^k$ (some of them $0$) such that, 
for all $1 \le i,j \le n$,
\[
(A_{ij} \o \pi_k)(x) = x^\al A_{ij}^k(x),
\]
and $A_{ij}^k(x) = B_{ij}^k(x) \ne 0$ for some $i,j$ and all $x \in W_k$.
By \eqref{ch}, the characteristic
polynomial of the $\cC$-family of normal matrices 
$A^k(x) = (A_{ij}^k(x))_{1 \le i,j \le n}$ has the form
\[
\ch(A^k(x))(\la) 
= (-1)^n \Big(\la^n + \sum_{j=2}^n (-1)^j x^{-j\al} a_j(\pi_k(x)) \la^{n-j}\Big).
\]
By theorem \ref{cCperturb},
there exists a finite covering $\{\pi_{kl} : W_{kl} \to W_k\}$ of $W_k$ (possibly shrinking $W_k$)
of the type described in $(2)$
such that, for all $l$, the family of polynomials 
$\ch(A^k \o \pi_{kl})$ admits a $\cC$-parameterization of its roots 
(the eigenvalues of $A^k \o \pi_{kl}$).
Eigenvectors of $(A^k \o \pi_{kl})(x)$ are also eigenvectors of $(A \o \pi_k \o \pi_{kl})(x)$ (and vice versa).

Let $l$ be fixed and let $y_0 \in W_{kl}$.
As there exist indices $1 \le i,j \le n$ such that $A_{ij}^k(x) \ne 0$ for all $x \in W_k$, and, thus,
$(A_{ij}^k \o \pi_{kl})(y_0) \ne 0$, not all eigenvalues of $(A^k \o \pi_{kl})(y_0)$ coincide.
Hence, $|\on{S}(\ch(A^k \o \pi_{kl})(y_0))| < |\on{S}(\ch(A(0)))|$, and 
the induction hypothesis implies the claim.
\end{demo}

The statement for Hermitian families $A$ can be proved in the same way, using theorem \ref{hypcCperturbth} 
instead of theorem \ref{cCperturb}.
\qed\end{demo}

\begin{remark}
The real analytic diagonalization of real analytic multiparameter families of symmetric matrices was treated by \cite[6.2]{KurdykaPaunescu08}. 
A one parameter version of theorem \ref{qamatrix} is proved in \cite{RainerAC} for 
$C^\infty$-curves of normal matrices $A$ such that $\ch(A)$ is everywhere normally nonflat (see \ref{1dim}). 
\end{remark}

If the parameter space is one dimensional, we recover a $\cC$-version of Rellich's classical perturbation result \cite[Satz 1]{Rellich37I}:

\begin{corollary}
Let $I \subseteq \R$ be an open interval. Consider a curve of Hermitian complex 
matrices 
\[
A(x)=(A_{ij}(x))_{1 \le i,j \le n},
\]
where the entries $A_{ij}$ (for $1 \le i,j \le n$) belong to $\cC(I,\C)$. 
Then there exist global $\cC$-parameterizations of the eigenvalues and the eigenvectors of $A$ on $I$. 
\end{corollary}

\begin{demo}{Proof}
The global statement for the eigenvectors can be proved by the arguments in the end of \cite[7.6]{AKLM98}.
\qed\end{demo}

\begin{example}[A nonflatness condition is necessary] \label{flat}
The following simple example (due to Rellich \cite{Rellich37I}, see also \cite[II 5.3]{Kato76}) shows that the above theorem 
is false if no nonflatness condition (such as quasianalyticity or normal nonflatness) is required: The eigenvectors of the smooth Hermitian family 
\[
A(x) := e^{-\frac{1}{x^2}}
\begin{pmatrix}
\cos \tfrac{2}{x} & \sin \tfrac{2}{x} \\
\sin \tfrac{2}{x} & - \cos \tfrac{2}{x}
\end{pmatrix} 
\text{ for } x \in \R \setminus \{0\}, \text{ and }  A(0) := 0,
\]   
cannot be chosen continuously near $0$. 
\end{example}

\begin{example}[Normality of $A$ is necessary] \label{normal}
Neither can the condition that $A$ is normal be omitted: 
Any choice of eigenvectors of the real analytic family
\[
A(x) := \begin{pmatrix}
0 & 1 \\
x & 0
\end{pmatrix}
\text{ for } x \in \R
\]
has a pole at $0$. The two parameter family 
\[
A(x,y) := \begin{pmatrix}
0 & x^2 \\
y^2 & 0
\end{pmatrix}
\text{ for } x,y \in \R
\]
has the eigenvalues $\pm xy$. But its eigenvectors cannot be chosen continuously near $0$, even after applying blow-ups or power substitutions.
\end{example}

\begin{theorem}[Regularity of the eigenvalues and eigenvectors] \label{regmatrix}
Let $M$ be a $\cC$-manifold.
Consider a family of normal complex 
matrices 
\[
A(x)=(A_{ij}(x))_{1 \le i,j \le n}
\]
(acting on a complex vector space $V = \C^n$), where the entries $A_{ij}$ (for $1 \le i,j \le n$) 
belong to $\cC(M,\C)$. For any compact subset $K \subseteq M$ there exists a relatively compact neighborhood $W$ of $K$ and 
parameterizations of the eigenvalues $\la_i$ and the eigenvectors $v_i$ (for $1 \le i \le n$) of $A$ on $W$ such that 
for all $i$:
\begin{enumerate}
\item[(1)] $\la_i \in \cW^{\cC}(W,\C)$ and $v_i \in \cW^{\cC}(W,\C^n)$.
\end{enumerate}
If $M$ is a open subset of $\R^q$, then: 
\begin{enumerate}
\item[(2)] $\la_i \in SBV(W,\C)$ and $v_i \in SBV(W,\C^n)$.
\end{enumerate}
\end{theorem}

\begin{demo}{Proof} 
The assertions for the eigenvalues follow immediately from the theorems
\ref{Wloc} and \ref{BVth}. 
The statements for the eigenvectors can be deduced from theorem \ref{qamatrix} in an analogous way as theorem 
\ref{Wloc} and theorem \ref{BVth} are deduced from theorem \ref{cCperturb} (compare with section \ref{secW11} and section \ref{secBV}).
\qed\end{demo}

\subsection{Example}
Consider the Hermitian family
\[
A(x,y) := \begin{pmatrix}
x & i y \\
-i y & -x
\end{pmatrix}
\text{ for } x,y \in \R.
\]
Its eigenvalues $\pm \sqrt{x^2+y^2}$ are not differentiable at $0$ and its eigenvectors cannot be arranged continuously near $0$.
Blowing up the origin, we end up with a family of Hermitian matrices which admits real analytic eigenvalues and eigenvectors; in coordinates:
\[
A(x,xy) = x \begin{pmatrix}
1 & i y \\
-i y & -1
\end{pmatrix} 
\]
has eigenvalues $\pm x \sqrt{1+y^2}$ and eigenvectors
\[
\binom{-1 - \sqrt{1+y^2}}{i y} \text{ and } \binom{i y}{-1 - \sqrt{1+y^2}};
\]
likewise,
\[
A(xy,y) = y \begin{pmatrix}
x & i  \\
-i  & -x
\end{pmatrix} 
\]
has eigenvalues $\pm y \sqrt{1+x^2}$ and eigenvectors
\[
\binom{-x + \sqrt{1+x^2}}{i} \text{ and } \binom{i}{-x + \sqrt{1+x^2}}.
\]
Setting
\begin{gather*}
v_1(x,y) := \binom{-1 - \sqrt{1+(\tfrac{y}{x})^2}}{i \tfrac{y}{x}}, ~ %\text{ and } 
v_2(x,y) := \binom{i \frac{y}{x}}{-1 - \sqrt{1+(\frac{y}{x})^2}}, ~  \text{ if } 0 < |y| \le |x|,\\
v_1(x,y) := \binom{- \tfrac{x}{y} + \sqrt{1+(\tfrac{x}{y})^2}}{i }, ~ %\text{ and } 
v_2(x,y) := \binom{i }{-\tfrac{x}{y} + \sqrt{1+(\tfrac{x}{y})^2}}, ~ \text{ if } 0 < |x| < |y|,\\ 
v_1(x,y) := \binom{1}{0}, ~ %\text{ and } 
v_2(x,y) := \binom{0}{1}, ~ \text{ if } y=0,\\
v_1(x,y) := \binom{1}{i}, ~ %\text{ and } 
v_2(x,y) := \binom{i}{1}, ~ \text{ if } x=0\ne y,
\end{gather*}
provides a choice of eigenvectors $v_1,v_2$ of $A$ which, clearly, is not continuous, but belongs to $\cW^{\cC}_{\on{loc}}$ 
(for any $\cC$ satisfying (\ref{cC}.1)--(\ref{cC}.6)) and, thus, also to $SBV_{\on{loc}}$.

\section{Applications to subanalytic functions} \label{secsuban}

\subsection{Subanalytic functions} 
Cf.\ \cite{BM88}.
Let $M$ be a real analytic manifold. A subset $X \subseteq M$ is called \emph{subanalytic} if each point of $M$ admits 
a neighborhood $U$ such that $X \cap U$ is a projection of a relatively compact semianalytic set.

Let $U$ be an open subanalytic subset of $\R^q$.  
Following \cite{Parusinski94} we call a function $f : U \to \R$ \emph{subanalytic} if the closure in $\R^q \times \R \mathbb P^1$ of the graph of $f$
is a subanalytic subset of $\R^q \times \R \mathbb P^1$.

Any continuous subanalytic function $f : U \to \R$ admits rectilinearization:
There exists a locally finite covering $\{\pi_k : U_k \to U\}$ of $U$, where each $\pi_k$ is a composite of finitely many mappings each of which is 
either a local blow-up with smooth center or a local power substitution, such that, for all $k$, the function $f \o \pi_k$ is real analytic \cite[1.4 \& 1.7]{BM90}.
This result was improved in \cite[2.7]{Parusinski94} to show that in the composition of the $\pi_k$ it is enough to substitute powers at the last step
after all local blow-ups.

\begin{theorem} \label{subanBV}
Let $U$ be an open subanalytic subset of $\R^q$.
Any continuous subanalytic function $f : U \to \R$ belongs to $\cW^{C^\om}_{\on{loc}}(U)$, and, thus, to $SBV_{\on{loc}}(U)$.
\end{theorem}

\begin{demo}{Proof}
This follows from rectilinearization and the reasoning in section \ref{secW11} and section \ref{secBV}.
\qed\end{demo}

\begin{theorem} \label{suban}
The roots of a family of polynomials $P$ whose coefficients are continuous subanalytic functions 
admit a parameterization in $\cW^{\cC^\om}_{\on{loc}}$,
and, thus, in $SBV_{\on{loc}}$.  
\end{theorem}

\begin{demo}{Proof}
Apply rectilinearization to the coefficients of $P$ and use theorem \ref{Comegath}. 
\qed\end{demo}

\begin{remark} \label{subanrk}
We cannot expect that for the rectilinearization of the roots of a continuous subanalytic hyperbolic family $P$ no local power substitutions are needed. 
This is shown by the following example:
\[
P(x)(z) := z^2 -|x|, \quad \text{for } x \in \R^q.
\]
If we additionally require that all coefficients of a subanalytic hyperbolic family $P$ are also arc-analytic, 
then indeed local blow-ups suffice,
by \cite[1.4]{BM90} (see also \cite[3.1]{Parusinski94}) and theorem \ref{hypcCperturbth}.
\end{remark}

\def\cprime{$'$}
\providecommand{\bysame}{\leavevmode\hbox to3em{\hrulefill}\thinspace}
\providecommand{\MR}{\relax\ifhmode\unskip\space\fi MR }
\providecommand{\MRhref}[2]{%
  \href{http://www.ams.org/mathscinet-getitem?mr=#1}{#2}
}
\providecommand{\href}[2]{#2}

%\bibliography{biblio}
%\bibliographystyle{amsplain}

\end{document}